\pdfoutput=1
\RequirePackage{ifpdf}
\ifpdf 
\documentclass[pdftex]{sigma}
\else
\documentclass{sigma}
\fi

\numberwithin{equation}{section}

\newtheorem{Theorem}{Theorem}[section]
\newtheorem*{Theorem*}{Theorem}
\newtheorem{Corollary}[Theorem]{Corollary}
\newtheorem{Lemma}[Theorem]{Lemma}
\newtheorem{Proposition}[Theorem]{Proposition}

\theoremstyle{definition}
\newtheorem{Definition}[Theorem]{Definition}

\newtheorem{Example}[Theorem]{Example}
\newtheorem{Remark}[Theorem]{Remark}

\newcommand{\SM}{(\mathbf{S}, \mathbf{M})}
\newcommand{\M}{\mathbf{M}}
\newcommand{\mathdash}{\relbar\mkern-9mu\relbar}

\begin{document}

\allowdisplaybreaks

\newcommand{\arXivNumber}{2312.06148}

\renewcommand{\PaperNumber}{019}

\FirstPageHeading

\ShortArticleName{Matrix Formulae and Skein Relations for Quasi-Cluster Algebras}

\ArticleName{Matrix Formulae and Skein Relations \\ for Quasi-Cluster Algebras}

\Author{Cody GILBERT~$^{\rm a}$, McCleary PHILBIN~$^{\rm b}$ and Kayla WRIGHT~$^{\rm c}$}

\AuthorNameForHeading{C.~Gilbert, M.~Philbin and K.~Wright}

\Address{$^{\rm a)}$~Department of Mathematics, Saint Louis University, Saint Louis, MO, USA}
\EmailD{\mail{cody.gilbert@slu.edu}}

\Address{$^{\rm b)}$~Department of Mathematics, University of Wisconsin - River Falls, River Falls, WI, USA}
\EmailD{\mail{mccleary.philbin@uwrf.edu}}

\Address{$^{\rm c)}$~Department of Mathematics, Johns Hopkins University, Baltimore, MD, USA}
\EmailD{\mail{kaylaw@jhu.edu}}

\ArticleDates{Received May 14, 2025, in final form January 30, 2026; Published online February 27, 2026}

\Abstract{In this paper, we give matrix formulae for non-orientable surfaces that provide the Laurent expansion for quasi-cluster variables, generalizing the orientable surface matrix formulae by Musiker--Williams. We additionally use our matrix formulas to prove the skein relations for the elements in the quasi-cluster algebra associated to curves on the non-orientable surface.}


\Classification{13F60; 05E16; 05C70}

\section{Introduction}

Cluster algebras were first defined by Fomin and Zelevinksy in the early 2000's in an effort to study problems regarding dual canonical bases and total positivity \cite{FZ02}. Since their axiomatization, endowing mathematical objects with a cluster structure has become a rapidly growing field of interest across many different disciplines in math and physics. For instance, Fomin and Shapiro defined a cluster structure associated to orientable topological surfaces (which actually generalized the deep geometric work of Fock and Goncharov in \cite{FG06, FG07} and Gekhtman, Shapiro, and Vainshtein \cite{GSV05}). The cluster algebra structure associated to an orientable surface has a completely topological description, but also can be bolstered to include more geometric significance. Namely, the cluster structure also arises from coordinate rings for the decorated Teichm\"uller space associated to the surface. The cluster variables can be thought of as hyperbolic lengths of geodesics on the surfaces also known as Penner coordinates \cite{P87}. Because of this interpretation, the Teichm\"uller theory was connected to combinatorics by associating elements in ${\rm PSL}_2(\mathbb{R})$ to arcs on the surface in order to obtain the Laurent expansions of the associated cluster variables. This result was first studied by Fock and Goncharov in \cite{FG06, FG07} in the coefficient-free case and generalized by Musiker and Williams in \cite{MW13}.

More recently, cluster-like structures were defined for non-orientable surfaces by Dupont and Palesi \cite{DP15}. They defined quasi-cluster algebras associated to non-orientable surfaces drawing inspiration from the orientable case in \cite{FST08}. The way they define their mutation or exchange relations is also inspired from the Teichm\"uller theory and geometry from \cite{FG06, FG07}. Positivity for quasi-cluster algebras was recently proven by Wilson through using snake and band graph combinatorics in \cite{Wil19} inspired by the original proof for positivity for cluster algebras from orientable surfaces in \cite{MSW11}.

In our paper, we create explicit matrix formulae for non-orientable surfaces that give the Laurent expansions for quasi-cluster variables in the associated quasi-cluster algebra. That is, we associate a product of matrices in ${\rm PSL}_2(\mathbb{R})$ to a curve on a non-orientable surface so that the Laurent expansion for the associated element of the quasi-cluster algebra can be extracted from the matrix. Our method provides the first non-recursive way to compute Laurent expansions for quasi-cluster variables directly on the non-orientable surface. Our process modifies the construction by Musiker and Williams \cite{MW13} by associating an element of ${\rm PSL}_2(\mathbb{R})$ when traversing through the non-trivial topology of a non-orientable surface. We in turn prove our result by relying on the combinatorics of snake and band graphs and a coefficient system using the language of laminations and Shear coordinates for non-orientable surfaces developed by Wilson in \cite{Wil19}. We also use our results to prove the skein relations for non-orientable surfaces.

Our paper is organized as follows: in Section~\ref{section:quasiclusteralgebras}, we review the cluster-like structure of non-orientable surfaces and topological notions we use throughout the paper. Section~\ref{section:wilsonexpansion} introduces the combinatorics of snake and band graphs in order to state the expansion formula for quasi-cluster variables as in \cite{Wil19}. Section~\ref{section:coefficients} introduces the coefficient system via laminations that give arbitrary coefficients for quasi-cluster algebras. Section~\ref{section:matrixorientable} recalls the matrix formulae by Musiker and Williams for cluster algebras from orientable surfaces and also extends their results to the context of principal laminations. This paves the way to our main results Theorems \ref{theorem:quasigoodmatchingenum} and \ref{thm:matrixformula} that is then proven in Section~\ref{section:expansionformulanonorientable}. Finally, Section~\ref{section:skeinrelations} proves the skein relations for non-orientable surfaces.

\section{Quasi-cluster algebras}\label{section:quasiclusteralgebras}
In this section, we review the definition of quasi-cluster algebras defined in \cite{DP15}. We give explicit computations for the so-called quasi-mutation rules that come from \textit{skein relations} as this was missing from the literature.
\begin{Definition} \label{definition::markedsurface}
Let $\mathbf{S}$ be a compact, connected smooth surface with boundary $\partial \mathbf{S}$. Let $\mathbf{M}$ be a finite set of points, we call \textit{marked points}, contained in $\mathbf{S}$ such that each connected component of $\partial \mathbf{S}$ has at least one point of $\mathbf{M}$. We say the pair $\SM$ is a \textit{marked surface}.
\end{Definition}

We call marked points on the interior of $\mathbf{S}$ \textit{punctures}. For the sake of clarity, in this article, we only deal with unpunctured surfaces, meaning marked points will only appear on the boundary.

\begin{Definition}\label{definition::regulararc}
A regular \textit{arc} $\tau$ in a marked surface $\SM$ is a curve in $\mathbf{S}$, considered up to isotopy relative its endpoints such that
\begin{enumerate}\itemsep=0pt
 \item[(1)] the endpoints of $\tau$ are $\mathbf{M}$,
 \item[(2)] $\tau$ has no self-intersections, except possibly at its endpoints,
 \item[(3)] except for its endpoints, $\tau$ does not intersect $\mathbf{M} \cup \partial \mathbf{S}$;
 \item[(4)] and $\tau$ does not cut out a monogon or bigon.
\end{enumerate}
We say the arc is a \textit{generalized arc} if $\tau$ intersects itself, dropping condition (2). We say the arcs that connect two marked points and lie completely on $\partial \mathbf{S}$ are \textit{boundary arcs}.
\end{Definition}

In our paper, we will be analyzing the specific case when $\mathbf{S}$ is non-orientable. By the classification of compact surface, any non-orientable surface is homeomorphic to the connected sum of $k$ projective planes. By this fact, we refer to $k$ as the \textit{non-orientable genus} of the surface. Recalling that the projective plane is a topological quotient of the 2-sphere by the antipodal map, we visualize these surfaces with \textit{crosscaps} $\bigotimes$. This symbol denotes the removal of a closed disk with the antipodal points identified.

On a non-orientable surface $\mathbf{S}$, a closed curve is said to be \textit{two-sided} if it admits an orientable regular neighborhood. If a closed curve does not admit such a neighborhood, it is said to be \textit{one-sided}. Since one-sided curves reverse the local orientation, they may only be contained in a~non-orientable surface.

As in the orientable case, an \textit{arc} is an isotopy class of a simple curve in $(\mathbf{S},\mathbf{M})$.

\begin{Definition}\label{definition:quasiarc}
 A \textit{quasi-arc} in $(\mathbf{S},\mathbf{M})$ is either an arc or a simple one-sided closed curve in the interior of $\mathbf{S}$.
\end{Definition}

We say that two arcs $\tau$, $\tau'$ in $\SM$ are \textit{compatible} if up to isotopy, they do not intersect one another. More formally, define $e(\tau, \tau')$ to be the minimum of the number of crossings of $\alpha$ and~$\alpha$' where~$\alpha$ is an arc isotopic to $\tau$ and~$\alpha'$ is an arc isotopic to $\tau'$
\[e(\tau, \tau') := \min\{\# \text{of crossings of } \alpha \text{ and } \alpha' \mid \alpha \simeq \tau \text{ and } \alpha' \simeq \tau'\},\]
\noindent where $\alpha$ ranges over arcs that are isotopic to $\tau$ and $\alpha'$ ranges over arcs isotopic to $\tau'$.

We say regular arcs $\tau$, $\tau'$ are \textit{compatible} if $e(\tau, \tau') = 0$. A maximal collection of pairwise compatible arcs in $\SM$ is called a \textit{quasi-triangulation of $\SM$}. If none of the arcs in this collection is a one-sided closed curve, we say it is simply a \textit{triangulation}.

All (quasi-)triangulations of a non-orientable surface are reachable via sequences of \textit{quasi-mutations} \cite[Proposition 3.12]{DP15}. These quasi-mutations are a larger class of local moves on the surface that are motivated by the \textit{skein relations}, studied in cluster algebra theory by \cite{FST08, MW13}. To discuss mutations, we require the notion of a quasi-seed.

Let $n$ be the \textit{rank}, i.e., the number of arcs in a (quasi)-triangulation and let $b$ the number of segments between marked points of the marked surface $(\mathbf{S}, \mathbf{M})$. Let $\mathcal{F}$ a field of rational functions in $n+b$ indeterminates. For each segment $a$ between marked points on $\partial \mathbf{S}$, we associate a variable $x_a\in \mathcal{F}$ such that $\{x_a\mid a\in\partial\mathbf{S}\}$ is algebraically independent in $\mathcal{F}$. We refer to~${\mathbb{Z}\mathbb{P}=\mathbb{Z}\big[x_a^{\pm} \mid a\in\partial\mathbf{S}\big\}\subset\mathcal{F}}$ as the \textit{ground ring}.

\begin{Definition}\label{definition:quasiseed}
 A quasi-seed associated with $(\mathbf{S}, \mathbf{M})$ in $\mathcal{F}$ is a pair $\Sigma=(T, \mathbf{x})$ such that
 \begin{itemize}\itemsep=0pt
 \item $T$ is a quasi-triangulation of $(\mathbf{S},\mathbf{M})$,
 \item $\mathbf{x}=\{x_t \mid t\in T\}$ is a free generating set of the field $\mathcal{F}$ over $\mathbb{Z}\mathbb{P}$.
 \end{itemize}
 The set $\{x_t \mid t\in T\}$ is called the \textit{quasi-cluster} of the quasi-seed $\Sigma$. We say a quasi-seed is a \textit{seed} if the corresponding quasi-triangulation is a triangulation.
\end{Definition}

\begin{Definition}\label{anti-self-folded-triangle}
 An \textit{anti-self-folded triangle} is any triangle of a quasi-triangulation with two edges identified by an orientation-reversing isometry.
\end{Definition}

For example, in Figure~\ref{fig:mutation1-example}, the triangle with sides $(t,t,a)$ is an anti-self-folded triangle.

Before stating the quasi-mutation rules, we derive two of the mutations from the skein relations. The other, significantly longer, computation is given in the appendix.

\begin{Example}
 Let $t$ be an arc in the anti-self-folded triangles with sides $(t,t,a)$, and let $t'$ be a~one-sided curve in an annuli with boundary $a$. Figure~\ref{fig:mutation1-example} demonstrates how to resolve the crossing of $t$ and $t'$ using the Ptolemy relation, which yields the relation $x_tx_{t'}=x_a$. We explain how to push a loop through the crosscap in the appendix. With this we can derive the second and third quasi-mutations given in Definition~\ref{definition:quasimutation}.
 \begin{figure}[t]\centering
 \includegraphics[width=.45\textwidth]{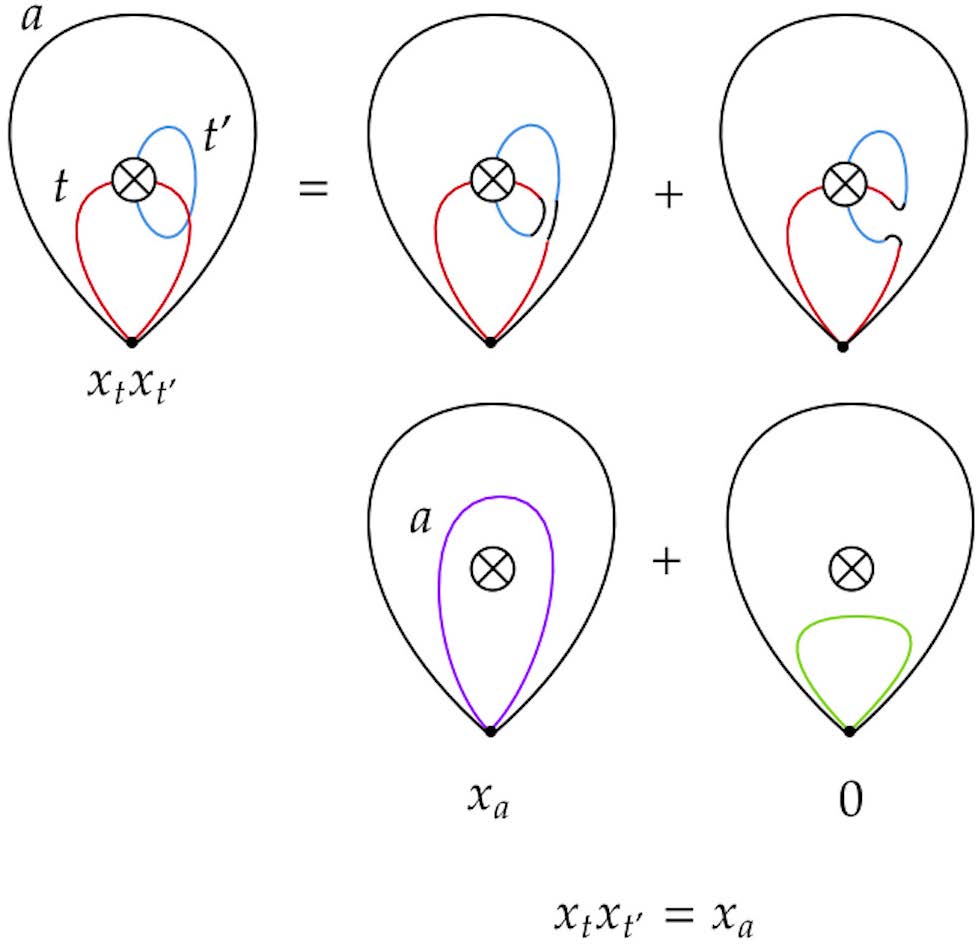}
 \caption{The skein relations used to define quasi-mutation (3).}
 \label{fig:mutation1-example}
 \end{figure}
 The justification for pushing the loop through the crosscap is given in the appendix.
\end{Example}

With this motivation, we have the following definitions of quasi-mutation.

\begin{Definition}\label{definition:quasimutation}
 Given $t\in T$, the \textit{quasi-mutation} of $\Sigma$ in the direction $t$ is the pair $\mu_t(T,\mathbf{x})=(T',\mathbf{x'})$, where $T'=\mu_t(T)=T\setminus \{t\}\sqcup\{t'\}$ and $\mathbf{x'}=\{x_v \mid v\in T'\}$ such that $x_{t'}$ is defined as follows:\looseness=-1
 \begin{enumerate}\itemsep=0pt
 \item[(1)] If $t$ is an arc separating two distinct triangles with sides $(a,b,t)$ and $(c,d,t)$, then the relation is given by the Ptolemy relation for arcs $x_tx_{t'}=x_ax_c+x_bx_d$,

 \smallskip

\centerline{\includegraphics[width=.37\textwidth]{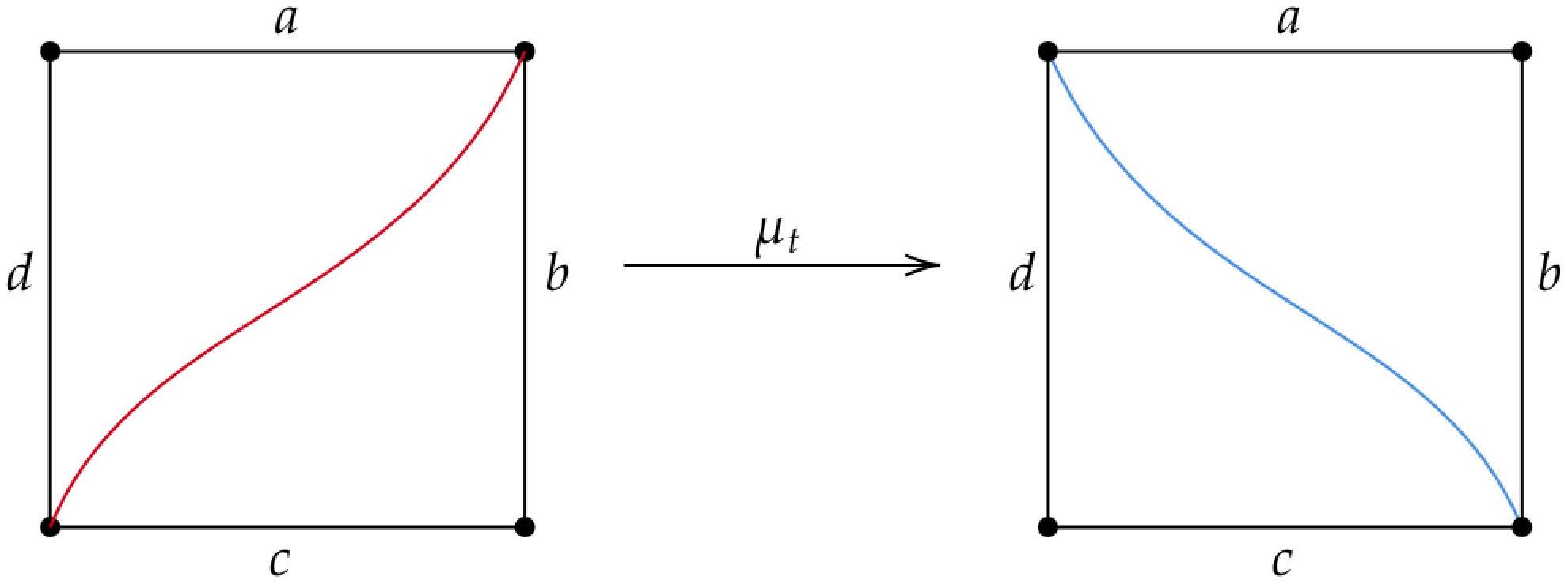}}

 \item[(2)] If $t$ is an arc in an anti-self-folded triangle with sides $(t,t,a)$, then the relation is $x_tx_{t'}=x_a$,

\centerline{\includegraphics[width=.37\textwidth]{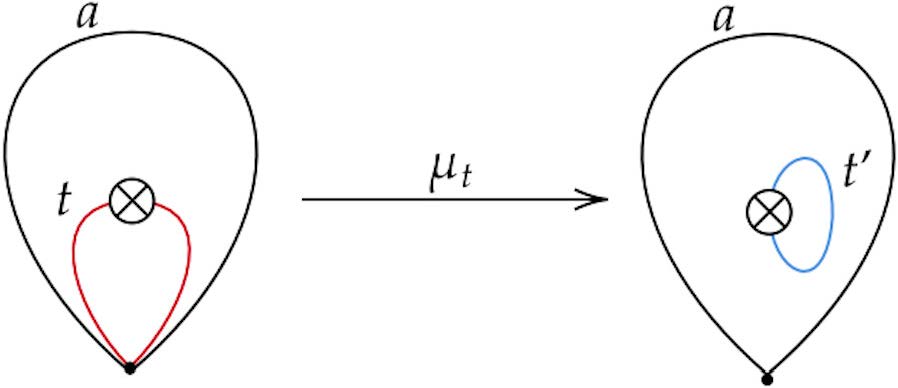}}

 \item[(3)] If $t$ is a one-sided curve in an annuli with boundary $a$, then the relation is $x_tx_{t'}=x_a$,

\centerline{\includegraphics[width=.37\textwidth]{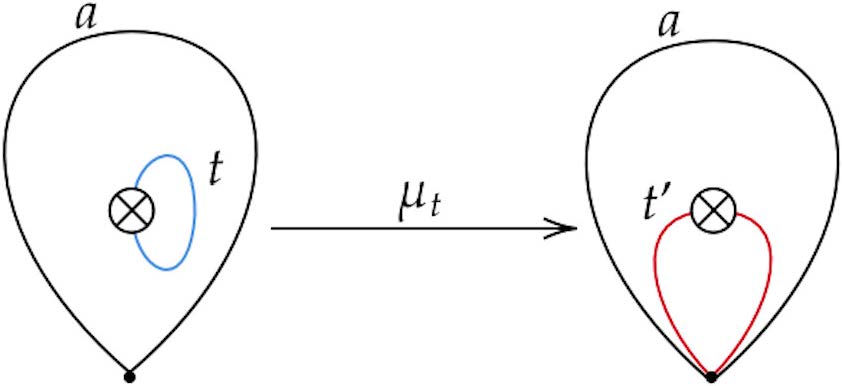}}

 \item[(4)] If $t$ is an arc separating a triangle with sides $(a,b,t)$ and an annuli with boundary $t$ and one-sided curve $d$, then the relation is $x_tx_{t'}=(x_a+x_b)^2+x_d^2x_ax_b$,

\centerline{\includegraphics[width=.47\textwidth]{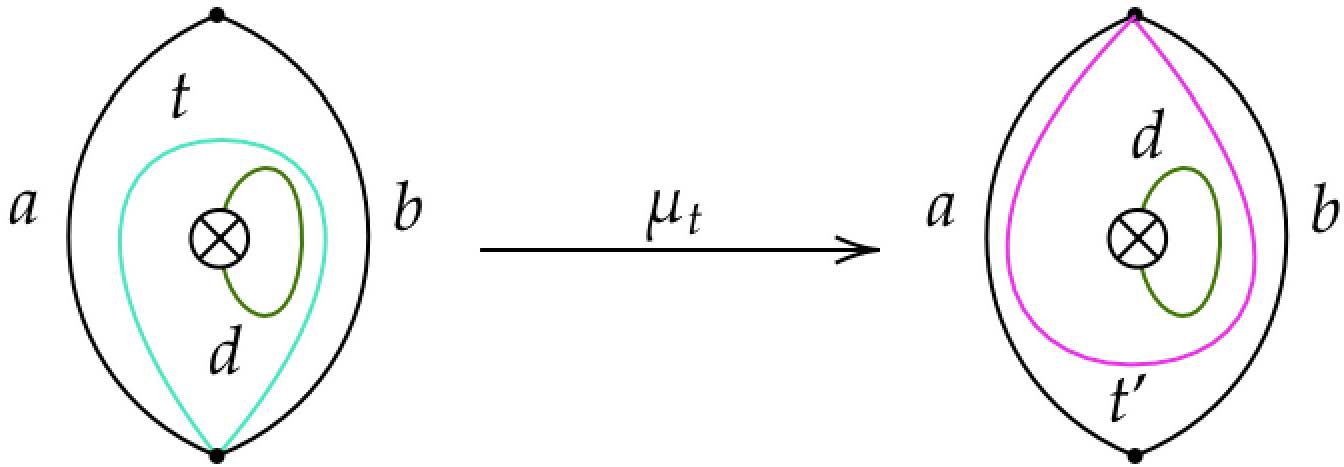}}
\end{enumerate}
\end{Definition}

Now that we have the necessary topological notions, we are ready to define the cluster structure for these non-orientable surfaces. This was first defined in \cite{DP15} which was inspired by work in \cite{FST08} in the orientable setting.

\begin{Definition}\label{definition:quasiclusteralgebra}
Let $\mathcal{X}$ be the collection of all quasi-cluster variables obtained by iterated quasi-mutation from an initial seed $\Sigma$. The \textit{quasi-cluster algebra} is the ring generated by the quasi-cluster variables over the ground ring $\mathbb{Z}\mathbb{P}$, i.e., $\mathcal{A}(\Sigma) = \mathbb{Z}\mathbb{P}[\mathcal{X}]$.
\end{Definition}

\section[Expansion formula via snake and band graphs for non-orientable surfaces]{Expansion formula via snake and band graphs\\ for non-orientable surfaces}\label{section:wilsonexpansion}

In this section, we review the expansion formulae for quasi-cluster variables in terms of snake and band graphs \cite{Wil19} without coefficients. We postpone the discussion on coefficients until Section~\ref{section:coefficients}. We begin by defining the notion of a snake graph as in~\cite{MSW11} and then review the definition of band graphs that come from non-orientable surfaces as in~\cite{Wil19}.
\begin{Definition}
 A \emph{tile} is a copy of the cycle graph $C_4$ on four vertices, embedded in $\mathbb{R}^2$ as a~square with four cardinal directions, see Figure~\ref{fig:tile}.
\begin{figure}[!ht]
 \centering
 \includegraphics[scale=.45]{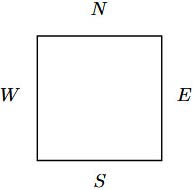}
 \caption{A tile with the four cardinal directions.}
 \label{fig:tile}
\end{figure}
\end{Definition}

We glue tiles together in a particular way to obtain a snake or band graph. In particular, a~snake or band graph can be thought of as a~sequence of tiles glued along either the North or East edge of the previous tile. We describe how to construct a~snake graph from an arc on a~triangulated orientable surface.

\begin{Definition}\label{definition:snakegraphofarc}
 Let $\gamma$ be an arc overlayed on a triangulation $T$ of an orientable surface $\SM$. The \emph{snake graph associated to $\gamma$} is a sequence of tiles $\mathcal{G}_\gamma = (G_1, \dots, G_d)$ where
 \begin{itemize}\itemsep=0pt
 \item the tiles $G_i$ correspond to the arc $\tau_i$ of $T$ that $\gamma$ intersects,
 \item and we attach $G_i$ to $G_{i+1}$ along the unique shared edge of the local quadrilaterals containing~$\tau_i$ and $\tau_{i+1}$.
 \end{itemize}
 Label the $N$, $E$, $S$, $W$ edges of each tile with the corresponding arcs in the local quadrilateral so that the relative orientation of the tiles alternate, see \cite{MSW11, Wil19} for details.
\end{Definition}

We now give the definition of a band graph associated to a one-sided closed curve as in \cite{Wil19}. We suppress explicit examples until later sections of the paper, see Figure~\ref{fig:MWBand}.

\begin{Definition}\label{definition:bandgraphquasiarc}
 Let $\gamma$ be a one-sided closed curve overlayed on a quasi-triangulation $T$, without quasi-arcs, of a non-orientable surface $\SM$. Let $x$ be a point on $\gamma \setminus T$. Consider the orientable double cover \smash{$\widetilde{\SM}$} of $\SM$ along the lifts of $T$, $\gamma$ and $x$. Note that $\tilde{T}$ is a triangulation of \smash{$\widetilde{\SM}$}, $\tilde{\gamma}$ is an orientable closed curve and $\tilde{x}_1$, $\tilde{x}_2$ of $x$ are antipodal points on $\tilde{\gamma}$.

 Let $\mathcal{G}_{\gamma, T, x} = (G_1, \dots, G_d)$ be the snake graph of the arc corresponding to tracing along $\tilde{\gamma}$ from~$\tilde{x}_1$ to $\tilde{x}_2$ clockwise. If we continue along $\tilde{\gamma}$ through one more intersection with $\tilde{T}$, we obtain~${(G_1, \dots, G_d, G_{d+1})}$, where $G_1$ and $G_{d+1}$ are lifts of the same local quadrilateral in $T$. Let~$b$ be the glued edge between $G_d$ and $G_{d+1}$ and $\overline{b}$ be the corresponding edge in $G_1$. The \emph{band graph associated to $\gamma$}, denoted $\mathcal{G}_{\gamma, x}$ is the snake graph $\mathcal{G}_{\tilde{\gamma}} = (G_1, \dots, G_d)$ with $b$ and $\overline{b}$ identified.
\end{Definition}

\begin{Remark}
 Up to isomorphism, the choice of point $x$ on the one-sided closed curve $\gamma$ does not affect the band graph or expansion formulae \cite{Wil19}.
\end{Remark}

To state the expansion formulae, we need the notion of (good) perfect matchings on snake (band) graphs.

\begin{Definition}\label{definition:goodperfectmatching}
 A \emph{perfect matching} of a snake graph $\mathcal{G}$ is a subset $P$ of the edges of $G$ that covers each vertex exactly once. A \emph{good perfect matching} $\overline{P}$ of a band graph $\overline{\mathcal{G}} = (G_1, \dots, G_d)$, glued along the edge $b$ from $x \mathdash y$, is a perfect matching $P$ of the un-identified underlying snake graph $\mathcal{G}$ such that one of the following holds:
 \begin{itemize}\itemsep=0pt
 \item $x$ and $y$ are matched to each other, i.e., $b \in \overline{P}$,
 \item $b \notin \overline{P}$ and $x$ and $y$ are matched with the two edges perpendicular to $b$ on tile $G_1$, or
 \item $b \notin \overline{P}$ and $x$ and $y$ are matched with the two edges perpendicular to $b$ on tile $G_d$.
 \end{itemize}
\end{Definition}

Now, we state the expansion formula for regular arcs without coefficients as in \cite[Theorem~5.3]{Wil19}.

\begin{Theorem}[\cite{Wil19}]\label{theorem:expformularegarcs}
 Let $\gamma$ be a regular arc overlayed on a quasi-triangulation $T$ without quasi-arcs of a non-orientable $\SM$. Let $\tilde{\gamma}$ be one of the two lifts of $\gamma$ on the orientable double cover~$\SM$. Let $\mathcal{G}_{\tilde{\gamma}}$ be the snake graph associated to $\tilde{\gamma}$ and $\mathcal{P}$ be the set of perfect matchings on $\mathcal{G}_{\tilde{\gamma}}$. Then the quasi-cluster variable $x_\gamma$ can be expressed as follows:
 \[x_\gamma = \frac{1}{\operatorname{cross}_T(\gamma)} \sum_{P \in \mathcal{P}} x(P),\]
 where cross$_T(\gamma) = x_{i_1} \cdots x_{i_d}$ is the crossing monomial which keeps track of the arcs of the triangulation $\tau_{i_1}, \dots, \tau_{i_d}$ that $\gamma$ crosses, counting multiplicities, and $x(P)$ is the product of all the edge labels of $\mathcal{G}_{\tilde{\gamma}}$ in $P$.
\end{Theorem}

The analogous expansion formula for one-sided closed curves without coefficients is given in~\cite[Theorem 5.19]{Wil19}.

\begin{Theorem}[{\cite{Wil19}}]\label{theorem:expformulaquasiarcsnocoeff}
 Let $\gamma$ be a one-sided closed curve overlayed on a quasi-triangulation $T$ without quasi-arcs of a non-orientable $\SM$. Let $\tilde{\gamma}$ be the lift of $\gamma$ on the orientable double cover $\SM$. Let $\mathcal{G}_{\tilde{\gamma}}$ be the band graph associated to $\tilde{\gamma}$ and $\mathcal{P}$ be the set of good perfect matchings on $\mathcal{G}_{\tilde{\gamma}}$. Then the quasi-cluster variable $x_\gamma$ can be expressed as follows:
 \[x_\gamma = \frac{1}{\operatorname{cross}_T(\gamma)} \sum_{P \in \mathcal{P}} x(P),\]
 where cross$_T(\gamma) = x_{i_1} \cdots x_{i_d}$ is the crossing monomial which keeps track of the arcs of the triangulation $\tau_{i_1}, \dots, \tau_{i_d}$ that $\gamma$ crosses, counting multiplicities, and $x(P)$ is the product of all the edge labels of $\mathcal{G}_{\tilde{\gamma}}$ in $P$.
\end{Theorem}

\section{Coefficients using principal laminations}\label{section:coefficients}

We review the arbitrary coefficients defined in \cite{Wil19} via principal laminations that complete the expansion formula mentioned in Section~\ref{section:wilsonexpansion}. We will then take inspiration from this coefficient system to define a poset structure of the set of good perfect matchings associated to a quasi-arc~$\gamma$ in a marked surface $\SM$.

\begin{Definition}\label{definition:lamination}
A set of self-non-intersecting and pairwise non-intersecting curves~$L$ on a~mar\-ked surface~$\SM$ is called a \emph{lamination} if each~$\ell \in L$ is any of the following curves:
\begin{itemize}\itemsep=0pt
 \item a one-sided closed curve,
 \item a two-sided closed curve that does not bound a disk, or a M\"obius strip,
 \item a curve that connected two points on $\partial \mathbf{S} \setminus \M$ that is not isotopic to a boundary arc.
\end{itemize}
A \emph{multilamination} $\mathbf{\mathcal{L}}$ of $\SM$ is a finite collection of laminations of $\SM$.
\end{Definition}

\begin{Example}
The figure below is an example of a lamination in red on the Klein bottle with four marked points.
\begin{figure}[!ht]
 \centering
 \includegraphics[scale=.45]{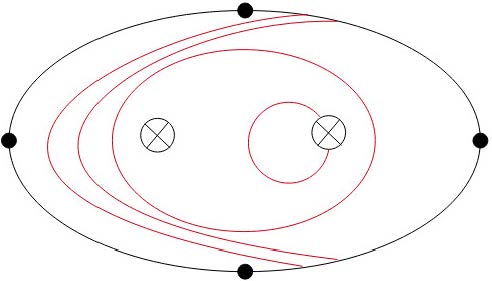}
\end{figure}
\end{Example}

\noindent We now define a principal lamination for marked surfaces in the unpunctured case.

\begin{Definition}\label{definition:principallaminations}
 Let $\gamma$ be an arc in $\SM$. We associate the lamination $L_\gamma$ to the arc $\gamma$ via the following:
 \begin{itemize}\itemsep=0pt
 \item If $\gamma$ is a regular arc, take $L_\gamma$ to be a lamination that runs along $\gamma$ in a small neighborhood, but turns clockwise (counterclockwise) at the marked point and ends when at the boundary.
 \item If $\gamma$ is a quasi-arc, take $L_\gamma$ to be either a lamination that runs along $\gamma$ in a small neighborhood and has endpoints on the boundary or take $L_\gamma$ to be the 1 sided closed curve that is compatible with $\gamma$.
 \end{itemize}
\end{Definition}

\begin{Remark}
 By definition, there are two choices for $L_{\gamma}$ for each $\gamma$.
\end{Remark}

\begin{Remark}
In the context of non-orientable surfaces, it is not entirely appropriate to describe direction via words like ``clockwise'' or ``counterclockwise''; however, by lifting to the orientable double cover, we can still retain these notions, for our purposes, as long as we are consistent with the Shear coordinates of $L_{\gamma}$ and $\gamma$, see Definition~\ref{definition:Shearcoordinate}.

In particular, the arc $\gamma$ on the non-orientable surface lifts to two arcs $\gamma'$ and $\gamma''$ on the orientable double cover, and these two arcs will have Shear coordinates with opposite sign.
\end{Remark}

\begin{figure}[!ht]
 \centering
\includegraphics[scale=.4]{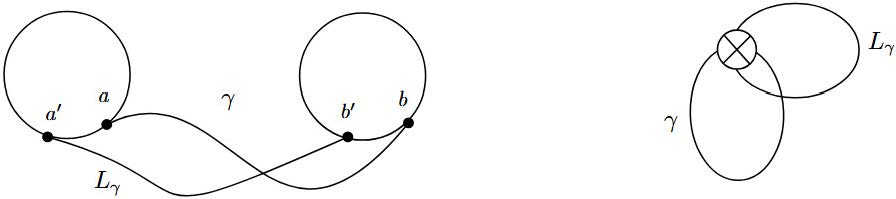}
\caption{Examples of the conditions for a laminations. On the left, the case where $\gamma$ is not a quasi-arc and on the right, the case where $\gamma$ is a quasi-arc.}\label{fig:laminationwitharc}
\end{figure}

\begin{Definition}\label{definition:principallamination}
 Let $T$ be a triangulation of $\SM$. Taking the collection of all laminations associated to the arcs $\tau \in T$ is called a \emph{principal lamination}. That is, a multilamination of the form $\mathcal{L}_T = \{L_\tau \mid \tau \in T\}$ is a principal lamination.
\end{Definition}

In order to define the coefficients seen in \cite{Wil19}, we must use the notion of a~principal lamination to define Shear coordinates, a~coordinate we place on the diagonal of a~local quadrilateral in a~triangulation.

\begin{figure}[t]
 \centering
 \includegraphics[scale=.44]{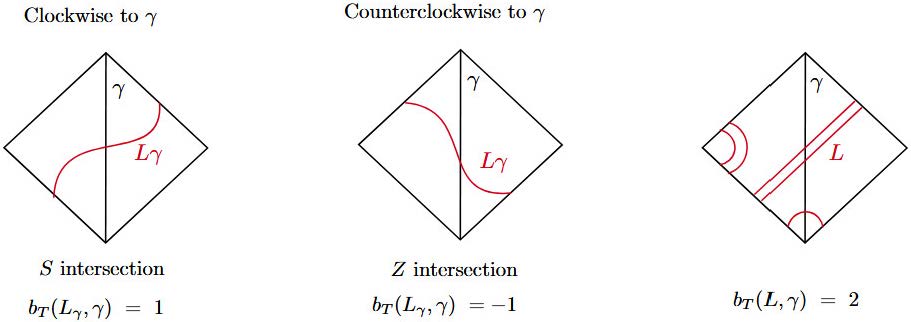}
 \caption{The definition of $S$ and $Z$ intersections on the leftmost and middle figures.} \label{fig:SZintersection}
\end{figure}

\begin{Definition}\label{definition:Shearcoordinate}
 Let $\gamma$ be an arc in some triangulation $T$ of $\SM$, let $L$ be a~lamination and let $Q_\gamma$ be the local quadrilateral that $\gamma$ is a diagonal of. The \emph{Shear coordinate of $L$ and $\gamma$ with respect to $T$}, denoted $b_T(L,\gamma)$, is given by
 \[b_T(L,\gamma) = \#\{S\text{-intersections of } L \text{ with } Q_\gamma\} - \#\{Z\text{-intersections of }L \text{ with }Q_\gamma\},\]
 where $S$-(resp.~$Z$-)intersections are illustrated in Figure~\ref{fig:SZintersection}.
\end{Definition}
We now use the notion of Shear coordinate to connect laminations to snake/band graphs and their perfect matchings. Namely, each diagonal or label of a~tile of a~snake/band graph corresponds to some arc in a triangulation. We use the following definition to assign a $\pm 1$ to an orientation of a diagonal associated to a tile in a snake/band graph.
\begin{Definition}
Let $\mathcal{G} = (G_1, \dots, G_d)$ be a snake (resp.\ band) graph and let $P$ be a~perfect matching (resp.\ good perfect matching) of $\mathcal{G}$. For each tile $G_i$ of $\mathcal{G}$ labeled on its diagonal, $P$ induces an \emph{orientation on the diagonal of $G_i$}. The orientation of $i$ is governed by the unique path from the SW vertex of $G_1$ to the NE vertex of $G_d$ taking alternating edges along $P$ and the diagonals of $\mathcal{G}$.
\end{Definition}

\begin{figure}[!ht]
 \centering
 \includegraphics[scale=.355]{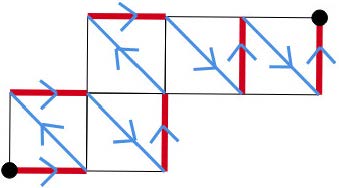}
 \caption{Orientation of diagonals induced by a perfect matching of a snake graph.}\label{fig:orientationofsnakegraphtiles}
\end{figure}

\begin{Example}
Consider the following perfect matching of the snake graph highlighted in red in Figure \ref{fig:orientationofsnakegraphtiles}. The orientation of each diagonal is designated by the path highlighted in blue. Note that this is also known in the literature as a $T$-path defined in \cite{S08}.
\end{Example}

The orientation of each tile in a snake graph allows for the following definition.

\begin{Definition}\label{definition:L_Toriented}
Let $\mathcal{G} = \mathcal{G}_{\alpha, T}$ be the snake (resp.\ band) graph associated to the curve~$\alpha$ and triangulation $T$ of $\SM$. Let $\mathcal{L}_T$ be a principal lamination of~$\SM$. A diagonal $\gamma_{i_k}$ of $\mathcal{G}$ is \emph{$\mathcal{L}_T$-oriented} with respect to a perfect matching (resp.\ good perfect matching)~$P$~if
\begin{itemize}\itemsep=0pt
 \item The tile is indexed odd and either
 \begin{itemize}\itemsep=0pt
 \item $b_T(L_{\gamma_{i_k}}, \gamma_{i_k}) = 1$ and the diagonal on the tile is oriented down, or
 \item $b_T(L_{\gamma_{i_k}}, \gamma_{i_k}) = -1$ and the diagonal on the tile is oriented up.
 \end{itemize}
 \item The tile is indexed even and either
 \begin{itemize}\itemsep=0pt
 \item $b_T(L_{\gamma_{i_k}}, \gamma_{i_k}) = 1$ and the diagonal on the tile is oriented up, or
 \item $b_T(L_{\gamma_{i_k}}, \gamma_{i_k}) = -1$ and the diagonal on the tile is oriented down.
 \end{itemize}
\end{itemize}
\end{Definition}

This definition of $\mathcal{L}_T$-oriented tells us exactly how to assign coefficients in our expansions.
\begin{Definition}\label{definition:wilsoncoefficients}
Given a (good) perfect matching $P$ of $\mathcal{G}$, the \emph{coefficient monomial} is given by
\[y_{\mathcal{L}_T}(P) = \prod_{\gamma_{i_k} \text{ is } \mathcal{L}_T \text{ oriented }} y_{\gamma_{i_k}}.\]
\end{Definition}

With these coefficients, we state Wilson's complete expansion formula \cite[Theorem 5.44]{Wil19} following the setup from Theorem~\ref{theorem:expformulaquasiarcsnocoeff} for quasi-arcs.

\begin{Theorem}[\cite{Wil19}]
 Let $\gamma$ be a one-sided closed curve overlayed on a quasi-triangulation $T$, let $\mathcal{L}_T$ be a principal lamination and $\mathcal{P}$ be the set of good perfect matchings of band graph $\mathcal{G}_{\overline{\gamma}}$. The Laurent expansion for the quasi-cluster variable $x_\gamma$ can be expressed as follows:
 \[x_\gamma = \frac{1}{\operatorname{bad}(\mathcal{L}_T,\gamma)}\frac{1}{\operatorname{cross}_T(\gamma)} \sum_{P \in \mathcal{P}} x(P) y_{\mathcal{L}_T}(P),\]
 where $x(P)$ is the product of all the edge labels of $\mathcal{G}_{\tilde{\gamma}}$ in $P$, $y_{\mathcal{L}_T}(P)$ is the coefficient monomial and bad$(\mathcal{L}_T,\gamma)$ is an error term counting the number of ``bad encounters'' as in {\rm\cite[\emph{Definition}~5.40]{Wil19}}.
\end{Theorem}

\begin{Remark}
 To say a bit more about the ``bad encounters'' mentioned in the previous theorem, suppose we have an elementary lamination $L_{\beta}$ that is a one-sided closed curve. If we consider the two lifts $\beta_1$ and $\beta_2$ of $\beta$ to the orientable double cover, we can then define the unique ``quasi-lamination'', $L^*_{\beta}$, associated to $L_{\beta}$ that lifts to two laminations that have the same shear coordinate. While $L^*_{\beta}$ is self-intersecting, and thus not a lamination, its two lifts are laminations of the orientable double cover.
 Given a principal lamination $\mathcal{L}_T$, we can define its quasi-principal lamination
 \begin{align*}
 \mathcal{L}_T^*={}&\mathcal{L}\setminus \{ L_{\beta} \mid L_{\beta}\in\mathcal{L}_{T} \text{ is a one-sided closed curve}\}\\
 &\cup \{ L^*_{\beta} \mid L_{\beta}\in\mathcal{L}_{T} \text{ is a one-sided closed curve}\}.
 \end{align*}
For any directed quasi-arc $\gamma$, we have $\gamma$ has a bad encounter with $\mathcal{L}^*_{\beta}$ if either $\gamma$ is a one-sided closed curve homotopic to $L_{\beta}$ or, within a small neighborhood, $\gamma$ intersects $\beta$, $L^*_{\beta}$, $L^*_{\beta}$, $\beta$, in that exact order. We then define
 \[\operatorname{bad}(\mathcal{\textbf{L}}_T^*, \gamma):=\prod_{\beta\in T} y_{\beta}^{a_{\beta}},\] where $a_b$ is the number of bad encounters of $\gamma$ with $L^*_{\beta}\in \mathcal{L}^*_T$. For more detail, read \cite[Section~5.3.2]{Wil19}.
\end{Remark}

\section[Matrix formulae on orientable surfaces with respect to principal laminations]{Matrix formulae on orientable surfaces\\ with respect to principal laminations}\label{section:matrixorientable}

In this section, we present our main results in Theorems \ref{thm:matrixformulalamination} and \ref{theorem:perfectmatchingenum}, and then our main results in Theorems \ref{theorem:quasigoodmatchingenum} and \ref{thm:matrixformula}.

We generalize the matrix formulae from \cite{MW13}, so they work when using coefficients from an arbitrary principal lamination, rather than just principal coefficients, that is, the case where ${b_T(\mathcal{L}_T, \gamma)=1}$ for all $\gamma\in T$. This will be necessary when we eventually generalize these matrix formulae to non-orientable surfaces, for if we lift an arc $\alpha$ on a non-orientable surface to an orientable surface, it will lift to two arcs $\alpha_1$ and $\alpha_2$, one of which will have an $S$-intersection, while the other will have a $Z$-intersection. With this in mind, many of the results on cluster algebras found in this section will naturally carry over to quasi-cluster algebras due to the existence of the orientable double cover.

Fix a marked surface $\SM$, triangulation $T$ of $\SM$ and principal lamination $\mathcal{L}_T=\{L_{\gamma} \mid \gamma\in T\}$. To a generalized arc or closed loop $\gamma$, we will be associating a cluster algebra element $x_{\gamma}$ using products of matrices in ${\rm PSL}_2(\mathbb{R})$. This can be achieved by recreating $\gamma$ via a~concatenation of various elementary steps, each of which has an associated matrix. The ensuing path created by adjoining these elementary steps will be referred to as an $M$-path. Upon taking the product of these matrices, $x_{\gamma}$ is either the upper right entry or the trace of the resulting matrix, depending on whether $\gamma$ was a generalized curve or a closed loop, respectively. The expansion for $x_{\gamma}$ achieved in this manner coincides with what one finds using more traditional means, such as mutation or the snake graph poset structure.

The definitions and results that we will be stating, and in some instances adjusting, from this section can be found in \cite[Sections~4 and~5]{MW13}. In order to define the aforementioned elementary steps and their matrices, we first need to establish some preliminary notation.

Elementary steps do not go from marked point to marked point, rather, they go between points which are close to marked points, but are not marked points themselves. With the above in mind, for each marked point $m\in \mathbf{M}$, we draw a small horocycle $h_m$ locally around~$m$, and as $m$ is on the boundary, we just consider $h_m\cap S$. We may assume the circles are small enough so that $h_m\cap h_{m'}=\varnothing$ for distinct marked points $m$, $m'$. For each arc $\tau\in T$, and marked point~${m\in \mathbf{M}}$ incident to $\tau$, we let~$v_{m, \tau}$ denote the intersection point $h_m\cap \tau$ and we let~$v^{+}_{m, \tau}$ (resp.\ $v^{-}_{m, \tau}$) denote a point on $h_m$ which is very close to $v_{m, \tau}$ but in the clockwise (resp.\ counterclockwise) direction from $v_{m, \tau}$.

We now give the definition of elementary steps from \cite{MW13} with a slight modification.

\begin{Definition}[elementary steps]\label{definition:elementarysteps} \quad
\begin{itemize}\itemsep=0pt
\item For the first type of step, we consider two arcs $\tau$ and $\tau'$ from $T$ which are both incident to a marked point $m$ and which form a triangle with third side $\sigma\in T$. Then the first step is a curve that goes between $\tau$ and $\tau'$ along $h_m$. The matrix associated to this step is
\[
\left[\begin{matrix}
1 & 0\\
\pm\dfrac{x_{\sigma}}{x_{\tau}x_{\tau'}} & 1
\end{matrix}\right],
\]
where the sign of $\frac{x_{\sigma}}{x_{\tau}x_{\tau'}}$ is positive if the orientation is clockwise and it is negative otherwise.

\item For the second type of step, we cross $\tau$ by following $h_m$ between $v^{+}_{m, \tau}$ and $v^{-}_{m, \tau}$. The associated matrix is
\[
\begin{bmatrix}
y_{\tau}^{\delta_{-1, b_T(\mathcal{L}_T, \tau)}} & 0\\
0 & y_{\tau}^{\delta_{1, b_T(\mathcal{L}_T, \tau)}}
\end{bmatrix}
\] if we travel clockwise and{\samepage
\[
\begin{bmatrix}
y_{\tau}^{\delta_{1, b_T(\mathcal{L}_T, \tau)}} & 0\\
0 & y_{\tau}^{\delta_{-1, b_T(\mathcal{L}_T, \tau)}}
\end{bmatrix}
\] otherwise. Here, $\delta$ is the Kronecker delta.}

\item For the third type of step, we travel along a path parallel to a fixed arc $\tau$ connecting two points $v^{\pm}_{m, \tau}$ and $v^{\mp}_{m', \tau}$ associated to distinct marked points $m$, $m'$. The associated matrix~is
\[\begin{bmatrix}
0 & \pm x_{\tau}\\
\mp \dfrac{1}{x_{\tau}} & 0
\end{bmatrix},\] where we use $x_r$, $\frac{-1}{x_r}$ if this step sees $\tau$ on the right and $-x_r$, $\frac{1}{x_r}$ otherwise.
\end{itemize}

\begin{Remark}
The first and third elementary steps are identical to those in \cite{MW13}. However, the second step is modified in order to reflect the underlying coefficient system using principal laminations.
\end{Remark}

With the adjusted type two step, we must formally restate and reprove the major results from~\cite[Section 5]{MW13} with the additional $b_T(\mathcal{L}_T, \tau)=-1$ case in mind, as there are a few differences that manifest.

\end{Definition}

By concatenating these elementary step segments to make a path, we can associate a matrix to an arc or a loop in the following way.

\begin{Definition}[{\cite{MW13}}]
Given $\SM, T$, and a generalized arc $\gamma$ in $\SM$ from $s$ to $t$, we choose a curve $\rho_{\gamma}$ satisfying the following:
\begin{itemize}\itemsep=0pt
\item It begins at points of the form $v^{\pm}_{s, \tau}$ and ends at points of the form $v^{\pm}_{t, \tau'}$, where $\tau$, $\tau'$ are arcs of $T$ incident to $s$ and $t$, respectively.
\item It is a concatenation of the elementary steps from Definition~\ref{definition:elementarysteps}, and is isotopic to the segment of $\gamma$ between $h_s\cap \gamma$ and $h_t\cap \gamma$.
\item The intersections of $\rho_{\gamma}$ with $T$ are in bijection with the intersections of $\gamma$ with $T$.
\end{itemize}
An analogous definition holds when $\gamma$ is a closed loop as well, with the exception that we must have that $\rho_{\gamma}$ is isotopic to $\gamma$. In either instance, we refer to $\rho_{\gamma}$ as the \textit{$M$-path}. If we have a decomposition $\rho_{\gamma}=\rho_t\circ\cdots \rho_1$ into elementary steps $\rho_i$, then we define $M(\rho_{\gamma})=M(\rho_{t})\cdots M(\rho_{2})M(\rho_{1})$, and the identity matrix is the matrix associated to the empty path.
\end{Definition}

\noindent \textbf{Notation:}
Given a $2\times 2$ matrix $M=(m_{ij})$, let $\operatorname{ur}(M)$ denote~$m_{12}$ (the upper-right entry of~$M$) and $\operatorname{tr}(M)$ denote the trace of $M$.

\begin{Remark}[\cite{MW13}]
The matrices corresponding to elementary steps of type two are not in $\text{SL}_2(\mathbb{R})$. For some applications, we can instead use the matrices
\[\begin{bmatrix}
 y_{\tau}^{-\frac{1}{2}} && 0\\
 0 && y_{\tau}^{\frac{1}{2}}
\end{bmatrix}\qquad\text{and}\qquad
\begin{bmatrix}
 y_{\tau}^{\frac{1}{2}} && 0\\
 0 && y_{\tau}^{-\frac{1}{2}}
\end{bmatrix}.\]
If we make this replacement, we write $\overline{M}(\rho_{\gamma})=\overline{M}(\rho_{t})\cdots \overline{M}(\rho_{2})\overline{M}(\rho_{1})$, which is a product of matrices in $\text{SL}_2(\mathbb{R})$.
\end{Remark}

While $M(\rho_{\gamma})$ and $\overline{M}(\rho_{\gamma})$ depend on their choice of $\rho_{\gamma}$, it turns out that their trace and upper right entry depend only on $\gamma$.

\begin{Lemma}[{\cite[Lemma 4.8]{MW13}}]\label{lemma: independenceofpath}
Let $\gamma_1$ and $\gamma_2$ be a generalized arc and a closed loop, respectively, without contractible kinks. Then if $\rho$ and $\rho'$ are two $M$-paths associated to $\gamma_1$, then $|{\operatorname{ur}(M((\rho))}|=|{\operatorname{ur}(M(\rho'))}|$. Analogously, for any two $M$-paths $\rho$, $\rho'$ associated to $\gamma_2$, we have $|{\operatorname{tr}(M(\rho))}|=|{\operatorname{tr}(M(\rho'))}|$.
\end{Lemma}

\begin{figure}[!ht]
 \centering
 \includegraphics[scale=.2]{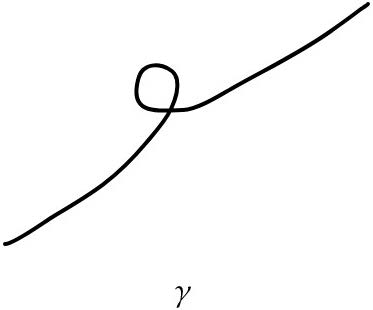}
 \caption{An example of a curve $\gamma$ with a contractible kink.}
\end{figure}

The lemma above allows us to make the following definition.

\begin{Definition}\label{definition:chilamination}
 Let $\gamma$ be a generalized arc and $\gamma'$ be closed loop, and let $\rho$ and $\rho'$ denote arbitrary $M$-paths associated to $\gamma$ and $\gamma'$, respectively. We associate the following algebraic quantities to $\gamma$ and $\gamma'$:
 \begin{enumerate}\itemsep=0pt
\item[(1)] $\chi_{\gamma, \mathcal{L}_T}=|{\operatorname{ur}(M(\rho))}|$ and $\chi_{\gamma', \mathcal{L}_T}=|{\operatorname{tr}(M(\rho'))}|$,
\item[(2)] $\overline{\chi}_{\gamma, \mathcal{L}_T}=|{\operatorname{ur}(\overline{M}(\rho))}|$ and $\overline{\chi}_{\gamma', \mathcal{L}_T}=|{\operatorname{tr}(\overline{M}(\rho'))}|$.
 \end{enumerate}
\end{Definition}

With the necessary definitions out of the way, we will now work towards proving the following result, which generalizes \cite[Theorem~4.11]{MW13} so that our matrix formulae now account for any choice of principal lamination $\mathcal{L}_T$.

\begin{Theorem}\label{thm:matrixformulalamination}
 Let $\SM$ be marked surface with a triangulation $T$ and let $T=\{\tau_1,\ldots, \tau_n\}$ be the corresponding triangulation. Let $\mathcal{A}_{\mathcal{L}_T}\SM$ be the corresponding cluster algebra associated to the principal lamination $\mathcal{L}_T$.
 \begin{itemize}\itemsep=0pt
\item Suppose $\gamma$ is a generalized arc in $S$ without contractible kinks. Let $\mathcal{G}_{T,\gamma}$ be the snake graph of $\gamma$ with respect to $T$. Then \[
 \chi_{\gamma, \mathcal{L}_T}=\frac{1}{\operatorname{cross}_T(\gamma)}\sum_{P}x(P)y_{\mathcal{L}_T}(P),
 \]
where the sum is over all perfect matchings $P$ of $\mathcal{G}_{T,\gamma}$. It follows that when $\gamma$ is an arc, then $\chi_{\gamma, \mathcal{L}_T}$ is the Laurent expansion of $x_\gamma$ with respect to $\mathcal{A}_{\mathcal{L}_T}\SM$.
 \item Suppose that $\gamma$ is a closed loop which is not contractible and has no contractible kinks. Then
\[\chi_{\gamma, \mathcal{L}_T}=\frac{1}{\operatorname{cross}_T(\gamma)}\sum_{P}x(P)y_{\mathcal{L}_T}(P),\] where the sum is over all good matchings $P$ of the band graph $\widetilde{\mathcal{G}}_{T, \gamma}$. Again, $\chi_{\gamma, \mathcal{L}_T}$ is the Laurent expansion of $\gamma$ with respect to $\mathcal{A}_{\mathcal{L}_T}\SM$.
 \end{itemize}
\end{Theorem}

To prove Theorem~\ref{thm:matrixformulalamination}, we will show $\chi_{\gamma, \mathcal{L}_T}$ coincides with the perfect (good) matching enumerator associated to their respective snake (band) graph. The process will be similar to that of~\cite{MW13}; however, with the additional $b_T(\gamma, \mathcal{L}_T)=-1$ case, we will need to make a few changes to the set-up as well as a few adjustments to the statements of the results. The proof of Theorem~\ref{thm:matrixformulalamination} can be found at the end of Section~\ref{section:StandardMPath}.

\subsection{Matchings of snake and band graphs}

To begin the proof, we will associate $2\times 2$ matrices to the parallelograms of a snake graph, which will then give us a way of representing the Laurent expansion of our graph in terms of products of matrices. The exposition defining the snake and band graphs is taken from \cite{MW13}.

\begin{Definition}[abstract snake graph]
An abstract snake graph with $d$ tiles is formed by concatenating the following pieces:
\begin{itemize}\itemsep=0pt
\item An initial triangle

\begin{figure}[!ht]
 \centering
 \includegraphics[scale=.5]{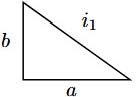}
\end{figure}

\item $d-1$ parallelograms $H_1, \ldots, H_{d-1}$, where each $H_j$ is either a north facing or east-pointing parallelogram

\begin{figure}[!ht]
 \centering
 \includegraphics[scale=.5]{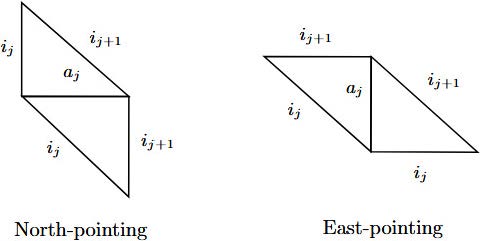}
\end{figure}

\item A final triangle based on whether $d$ is odd or even

\begin{figure}[!ht]
 \centering
 \includegraphics[scale=.5]{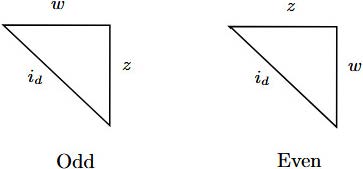}
\end{figure}

\end{itemize}
We then erase all diagonal edges from the figure.
\end{Definition}

\begin{Definition}[abstract band graph]\label{definition:bandgraph}
An abstract band graph for a one-sided curve with $d$ tiles is formed by concatenating the following puzzle pieces:
\begin{itemize}\itemsep=0pt
\item an initial triangle

\begin{figure}[!ht]
 \centering
 \includegraphics[scale=.5]{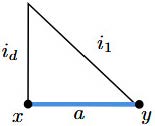}
\end{figure}

\item $d-1$ parallelograms $H_1, \ldots, H_{d-1}$, where each $H_j$ is as before,

\item a final triangle based on whether $d$ is odd or even.
\end{itemize}
\begin{figure}[!ht]
 \centering
 \includegraphics[scale=.5]{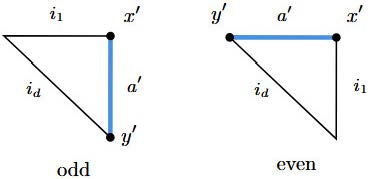}
\end{figure}

We then identify the edges $a$ and $a'$, the vertices $x$ and $x'$, and the vertices $y$ and $y'$. Lastly, we erase all diagonal edges from the figure.
\end{Definition}

Just like a traditional snake and band graph, we can consider perfect matchings of an abstract snake graph and good matchings of an abstract band graph. Furthermore, we can associate to each perfect/good matching its weight and coefficient monomials $x(P)$ and $y_{\mathcal{L}_{T}(P)}$ by letting~${T=\{\tau_{i_1},\ldots,\tau_{i_d}\}}$.

\begin{Definition}\label{definition:tilematrix}
Let $G$ be an abstract snake or band graph with $d$ tiles. For each parallelogram~$H_i$ of $G$, we associate a $2\times 2$ matrix $m_i$, where $m_i$ is either
\begin{gather*}
\begin{bmatrix}
1 & 0\\
\dfrac{x_{a_j}}{x_{i_j}x_{i_{j+1}}} & y_{i_j}
\end{bmatrix}\qquad\text{or}\quad
\begin{bmatrix}
\dfrac{x_{i_{j+1}}}{x_{i_j}} & x_{a_j}y_{i_j}\\
0 & \dfrac{x_{i_j}y_j}{x_{i_{j+1}}}
\end{bmatrix},\qquad\text{when}\quad b_{T}(\mathcal{L}_{T}, i_j)=1,
\\
\begin{bmatrix}
y_{i_j} & 0\vspace{1mm}\\
\dfrac{x_{a_j}y_{i_j}}{x_{i_j}x_{i_{j+1}}} & 1
\end{bmatrix}\qquad\text{or}\quad
\begin{bmatrix}
\dfrac{x_{i_{j+1}}y_{i_j}}{x_{i_j}} & x_{a_j}\\
0 & \dfrac{x_{i_j}}{x_{i_{j+1}}}
\end{bmatrix},\qquad\text{when}\quad b_{T}(\mathcal{L}_{T}, i_j)=-1.
\end{gather*}

In either case, the following conditions determine which of the two matrices we use:
\begin{itemize}\itemsep=0pt
\item $m_1$ is of the first type if $H_1$ is a north-pointing parallelogram, and otherwise it is of the second type,
\item for $i>1$, $m_i$ is of the first type if both $H_{i-1}$ and $H_i$ have the same shape, and otherwise, it is of the second type.
\end{itemize}
We can then define a matrix $M_d$ associated to $G$ defined by $M_d=m_{d-1}\cdots m_1$ when $d>1$. Otherwise, $M_1=I_2$.
\end{Definition}

\begin{Remark}
Definition~\ref{definition:tilematrix} differs from \cite{MW13} as it depends on the particular choice of Shear coordinate. This is necessary to reflect the coefficient system provided by the principal laminations in the non-orientable setting.
\end{Remark}

The upcoming sequence of results generalize Proposition 5.5, Corollary 5.6 and Theorem~5.4 from \cite{MW13} to principal laminations. Despite the additional case that stems from having $Z$-shape intersections, the proof techniques largely stay the same and the statements generalize in an expected manner.

\begin{Proposition}\label{prop:snakematrix}
Let $G$ be an abstract snake graph with $d$ tiles. Write \[M_d=\begin{bmatrix}
A_d & B_d\\
C_d & D_d
\end{bmatrix}.\]

For $d\geq 2$,{\samepage
\begin{itemize}\itemsep=0pt
\item when $b_{T}(\mathcal{L}_{T}, {i_{d-1}})=1$, we have
\begin{alignat*}{3}
& A_d = \frac{\sum_{P\in S_A}x(P)y_{\mathcal{L}_T}(P)}{(x_{i_1}\cdots x_{i_{d-1}})x_ax_w}, \qquad && B_d = \frac{\sum_{P\in S_B}x(P)y_{\mathcal{L}_T}(P)}{(x_{i_2}\cdots x_{i_{d-1}})x_bx_w},&\\
& C_d = \frac{\sum_{P\in S_C}x(P)y_{\mathcal{L}_T}(P)}{(x_{i}\cdots x_{i_d})x_ax_zy_{i_d}}, \qquad && D_d = \frac{\sum_{P\in S_D}x(P)y_{\mathcal{L}_T}(P)}{(x_{i_2}\cdots x_{i_{d}})x_bx_zy_{i_d}},&
\end{alignat*}
\item when $b_{T}(\mathcal{L}_{T}, i_{d-1})=-1$, we have
\begin{alignat*}{3}
& A_d = \frac{\sum_{P\in S_A}x(P)y_{\mathcal{L}_T}(P)}{(x_{i_1}\cdots x_{i_{d-1}})x_ax_wy_{i_d}} ,\qquad&& B_d = \frac{\sum_{P\in S_B}x(P)y_{\mathcal{L}_T}(P)}{(x_{i_2}\cdots x_{i_{d-1}})x_bx_w y_{i_d}},&\\
& C_d = \frac{\sum_{P\in S_C}x(P)y_{\mathcal{L}_T}(P)}{(x_{i}\cdots x_{i_d})x_ax_z} ,\qquad&& D_d = \frac{\sum_{P\in S_D}x(P)y_{\mathcal{L}_T}(P)}{(x_{i_2}\cdots x_{i_{d}})x_bx_z}.&
\end{alignat*}
\end{itemize}
For} $d=1$, the formulae remain the same, but the cases are determined by the sign of $b_T(\mathcal{L}_T, i_1)$.

Here, $S_A$, $S_B$, $S_C$ and $S_D$ are the sets of perfect matchings of $G$ which use the edges $\{a, w\}$, $\{b, w\}$, $\{a, z\}$ and $\{b, z\}$, respectively.
\end{Proposition}
\begin{proof}
 The case where $b_{T}(\mathcal{L}_{T}, {i_{d-1}})=1$ is exactly \cite[Proposition~5.5]{MW13}, so we only consider the case where $b_{T}(\mathcal{L}_{T}, {i_{d-1}})=-1$. Additionally, we note that differences between the two cases will only involve changes in $y_{\mathcal{L}_T}(P)$. This is due to the fact that, despite the change in principal lamination, the perfect matchings on the abstract snake graph remain the exact same, meaning~$x(P)$ will not change when flipping from a $S$-shape intersection to a $Z$-shape intersection, or vice versa.

 We perform a proof by induction on $d$, and consider what happens when one adds one more tile to a snake graph. When $d=1$, we just have a single tile
 meaning $S_A=\{a, w\}$, $S_B=S_C=\varnothing$ and $S_D=\{b,z\}$. When $b_T(\mathcal{L}_T, i_1)=-1$, we get
 \[
 A_1=\frac{x_ax_wy_{i_1}}{x_ax_wy_{i_1}} , \qquad B_1=0 , \qquad C_1=0, \qquad D_1=\frac{x_bx_z}{x_bx_z} ,
 \] hence $m_1=I_2$, as desired.
\begin{figure}[!ht]
 \centering
 \includegraphics[scale=.5]{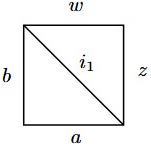}
\end{figure}

Just like \cite[Proposition 5.5]{MW13}, for $d>1$, we let $G'$ denote the graph obtained by gluing together the initial triangle, parallelograms $H_1,\ldots, H_{d-2}$ and the final triangle. For convenience, we~label the final triangle in $G'$ with $w'$ and $z'$ and note the orientation of this triangle depends on whether~${(d-1)}$ is odd or even. By changing labels, the graph $G'$ is isomorphic to the subgraph of $G$ consisting of the first~$d-1$ tiles. In particular, we either replace the edge label $w'$ with~$a_{d-1}$ and $z$ with $i_d$, or vice versa.

With this in mind, when considering $d=2$, we have may assume it has the following shape and labeling:
\begin{figure}[!ht]
 \centering
 \includegraphics[scale=.4]{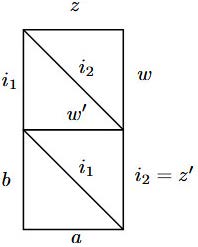}
\end{figure}

One can quickly check that when $b_T(\mathcal{L}_T, i_1)=1$ and $b_T(\mathcal{L}_T, i_2)=\pm 1$, then \[
M_2=
\begin{bmatrix}
 1 && 0\\
 \dfrac{x_{w'}}{x_{i_1}x_{i_2}} && y_{i_1}
\end{bmatrix}.
\]
Meanwhile, when $b_T(\mathcal{L}_T, i_1)=-1$ and $b_T(\mathcal{L}_T, i_2)=\pm 1$, then \[
 M_2=
\begin{bmatrix}
 y_{i_1} && 0\vspace{1mm}\\
 \dfrac{x_{w'}y_{i_1}}{x_{i_1}x_{i_2}} && 1
\end{bmatrix}.
\]
 Regardless, both matrices match the equations found in Definition~\ref{definition:tilematrix}.

To finish the proof, we simply observe that the bijections between the various perfect matchings utilized in the proof of \cite[Proposition 5.5]{MW13} apply here as well, regardless of our principal lamination. As such, by applying induction, we may conclude the following. Under the initial labeling scheme (used in the base cases), we obtain the following two cases depending on the sign of $b_T(\mathcal{L}_T, i_{d-1})$
\begin{align*}
 \begin{bmatrix}
 A_d & B_d\\
 C_d & D_d
 \end{bmatrix}&=
\begin{bmatrix}
 1 & 0\vspace{1mm}\\
 \dfrac{x_{a_{d-1}}}{x_{i_{d-1}}x_{i_d}} & y_{i_{d-1}}
\end{bmatrix}
\begin{bmatrix}
 A_{d-1} & B_{d-1}\\
 C_{d-1} & D_{d-1}
 \end{bmatrix},
 \\
 \begin{bmatrix}
 A_d & B_d\\
 C_d & D_d
 \end{bmatrix}&=
\begin{bmatrix}
 y_{i_{d-1}} & 0\vspace{1mm}\\
 \dfrac{x_{a_{d-1}}y_{i_{d-1}}}{x_{i_{d-1}}x_{i_d}} & 1
\end{bmatrix}
\begin{bmatrix}
 A_{d-1} & B_{d-1}\\
 C_{d-1} & D_{d-1}
 \end{bmatrix}.
\end{align*}

Likewise, under the other labeling scheme we obtain the following cases, again, depending on the sign of $b_T(\mathcal{L}_T, i_{d-1})$
\begin{align*}
 \begin{bmatrix}
 A_d & B_d\\
 C_d & D_d
 \end{bmatrix}&=
\begin{bmatrix}
 \dfrac{x_{i_d}}{x_{i_{d-1}}} & x_{a_{d-1}}y_{i_{d-1}}\vspace{1mm}\\
 0 & \dfrac{x_{i_{d-1}}y_{i_{d-1}}}{x_{i_d}}
\end{bmatrix}
\begin{bmatrix}
 A_{d-1} && B_{d-1}\\
 C_{d-1} && D_{d-1}
 \end{bmatrix},
 \\
 \begin{bmatrix}
 A_d & B_d\\
 C_d & D_d
 \end{bmatrix}&=
\begin{bmatrix}
 \dfrac{x_{i_d}y_{d-1}}{x_{i_{d-1}}} & x_{a_{d-1}}\vspace{1mm}\\
 0 & \dfrac{x_{i_{d-1}}}{x_{i_d}}
\end{bmatrix}
\begin{bmatrix}
 A_{d-1} & B_{d-1}\\
 C_{d-1} & D_{d-1}
 \end{bmatrix}.
\end{align*}
Comparing these two equations with Definition~\ref{definition:tilematrix}, we see the various cases agree with the cases of the definition of $M_d$.
\end{proof}

\begin{Corollary}
Let $G$ be an abstract snake graph with $d$ tiles. Write \smash{$M_d=\big[\begin{smallmatrix}
 A_d & B_d\\
C_d & D_d
\end{smallmatrix}\big]$}.
\begin{itemize}\itemsep=0pt
\item When $b_{T}(\mathcal{L}_{T}, i_{d})=1$, then
\begin{align*}
\frac{\sum_P x(P)y_{\mathcal{L}_T}(P)}{x_{i_1}\cdots x_{i_d}}&=\frac{x_a x_w A_d}{x_{i_d}}+\frac{x_b x_w B_d}{x_{i_1} x_{i_d}}+x_{a} x_z y_{i_d}C_d+\frac{x_b x_z y_{i_d}D_d}{x_{i_1}}\nonumber\\
&=\operatorname{ur}\left(\begin{bmatrix}
\dfrac{x_w}{x_{i_d}} & x_{z}y_{i_d}\vspace{1mm}\\
-\dfrac{1}{x_{z}} & 0
\end{bmatrix}
\begin{bmatrix}
A_d & B_d\\
C_d & D_d
\end{bmatrix}
\begin{bmatrix}
0 &x_a\\
-\dfrac{1}{x_a} & \dfrac{x_b}{x_{i_1}}
\end{bmatrix}\right).
\end{align*}

\item When $b_{T}(\mathcal{L}_{T}, i_{d})=-1$, then
\begin{align*}
\frac{\sum_P x(P)y_{\mathcal{L}_T}(P)}{x_{i_1}\cdots x_{i_d}}&=\frac{x_a x_w y_{i_d} A_d}{x_{i_d}}+\frac{x_b x_w y_{i_d} B_d}{x_{i_1} x_{i_d}}+x_{a} x_z C_d+\frac{x_b x_z D_d}{x_{i_1}}\\
&=\operatorname{ur}\left(\begin{bmatrix}
\dfrac{x_w y_{i_d}}{x_{i_d}} & x_{z}\vspace{1mm}\\
-\dfrac{y_{i_d}}{x_{z}} & 0
\end{bmatrix}
\begin{bmatrix}
A_d & B_d\\
C_d & D_d
\end{bmatrix}
\begin{bmatrix}
0 & x_a\\
-\dfrac{1}{x_a} & \dfrac{x_b}{x_{i_1}}
\end{bmatrix}\right), 
\end{align*}
where the sum is over all perfect matchings of $G$.
\end{itemize}
\end{Corollary}

The above corollary immediately implies the next theorem.

\begin{Theorem}\label{theorem:perfectmatchingenum}
 Suppose $G$ is an abstract snake graph with $d$ tiles. Then its perfect matching enumerator is given by
 \begin{itemize}\itemsep=0pt
 \item \[
 \sum_P x(P)y_{\mathcal{L}_T}(P)=x_{i_1}\cdots x_{i_d} \operatorname{ur}\left(\begin{bmatrix}
\dfrac{x_w}{x_{i_d}} & x_{z}y_{i_d}\vspace{1mm}\\
-\dfrac{1}{x_{z}} & 0
\end{bmatrix}
M_d
\begin{bmatrix}
0 & x_a\\
-\dfrac{1}{x_a} & \dfrac{x_b}{x_{i_1}}
\end{bmatrix}\right),\] when $ b_{T}(\mathcal{L}_{T}, i_{d})=1$.

\item \[
\sum_P x(P)y_{\mathcal{L}_T}(P)=x_{i_1}\cdots x_{i_d} \operatorname{ur}\left(\begin{bmatrix}
\dfrac{x_w y_{i_d}}{x_{i_d}} & x_{z}\vspace{1mm}\\
-\dfrac{y_{i_d}}{x_{z}} & 0
\end{bmatrix}
M_d
\begin{bmatrix}
0 & x_a\\
-\dfrac{1}{x_a} & \dfrac{x_b}{x_{i_1}}
\end{bmatrix}\right),\] when $b_{T}(\mathcal{L}_{T}, i_{d})=-1$.
\end{itemize}
\end{Theorem}

By adjusting our labeling by substituting $i_1$ for $w$, $a'$ for $z$, and $i_d$ for $b$ (see Definition~\ref{definition:bandgraph}), we obtain a sequence of statements and results for abstract band graphs that will be analogous to \cite[Proposition 5.7 and Corollary 5.8]{MW13}, but with respect to principal laminations. We end up with the following theorem.

\begin{Theorem}\label{theorem:goodmatchingenum}
 Suppose $G$ is an abstract band graph with $d$ tiles. Then its good matching enumerator, i.e., the weighted sum of good perfect matchings of $G$ with respect to $\mathcal{L}_T$ is given by
 \begin{itemize}\itemsep=0pt
 \item \[
 \sum_P x(P)y_{\mathcal{L}_T}(P)=x_{i_1}\cdots x_{i_d} \operatorname{tr}\left(\begin{bmatrix}
\dfrac{x_{i_1}}{x_{i_d}} & x_{a}y_{i_d}\\
0 & \dfrac{x_{i_d}y_{i_d}}{x_{i_1}}
\end{bmatrix}
M_d\right),\] when $ b_{T}(\mathcal{L}_{T}, i_{d})=1$.

\item \[
\sum_P x(P)y_{\mathcal{L}_T}(P)=x_{i_1}\cdots x_{i_d} \operatorname{tr}\left(\begin{bmatrix}
\dfrac{x_{i_1}y_{i_d}}{x_{i_d}} & x_{a}\\
0 & \dfrac{x_{i_d}}{x_{i_1}}
\end{bmatrix}
M_d\right),
\] when $ b_{T}(\mathcal{L}_{T}, i_{d})=-1$.
\end{itemize}
\end{Theorem}

\subsection[The standard M-path]{The standard $\boldsymbol{M}$-path}\label{section:StandardMPath}

Arcs do not have a unique associated $M$-path, but there is an algorithm for assigning a ``standard $M$-path'' $\rho_{\gamma}$ to a curve $\gamma$ that we will utilize in order to facilitate proofs. Here we recall the definition of a standard $M$-path for both an arc and a closed loop. Despite working with principal laminations for coefficients, the standard $M$-paths for generalized arcs and closed curves will not change. The only difference that will occur on the $M$-path side will involve adjusting our step two matrix, depending on whether the arc's corresponding lamination is an $S$-intersection or a $Z$-intersection.

\begin{figure}[!ht]
\centering
\includegraphics[scale=.55]{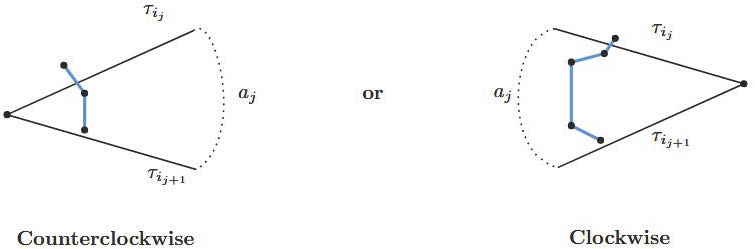}
\caption{Possible transitions between $\tau_{i_j}$ and $\tau_{i_{j+1}}$ in the standard $M$-path.}\label{fig:transitions}
\end{figure}

\begin{Definition}[standard $M$-paths]\label{definition:Mpath}
The definitions below are \cite[Definitions 5.9 and 5.11]{MW13}. In Definition~\ref{definition:Mpath2}, we generalize this construction for one-sided closed curves.

\textbf{Generalized arcs:} We begin with the case where we have a generalized arc $\gamma$. Let $\gamma$ be a~generalized arc that goes from point $P$ to a point $Q$, crossing the arcs $\tau_{i_1},\ldots,\tau_{i_d}$ in order. We label the initial triangle $\Delta_0$ that intersects $\gamma$ with sides $a$, $b$, and $\tau_{i_1}$ in clockwise order so that~$P$ is the intersection of $a$ and $b$. We label the final crossed triangle $\Delta_d$ with sides $w$, $z$ and $\tau_{i_d}$ in clockwise order, with $Q$ being the intersection of the arcs $w$ and $z$.

The path $\rho_{\gamma}$ begins at the point $v^{\pm}_{P,a}$, where the sign is chosen so the point lies in the triangle. The first two steps of $\rho_{\gamma}$ consists of traveling along $a$ and then jumping from $a$ to $\tau_{i_1}$, i.e., we perform a step three and then a step one move. After these two moves, we are at a point of the form $v^{\pm}_{*, \tau_{i_1}}$, but we have not yet crossed $\tau_{i_1}$.

Until we reach the final triangle, we follow the steps illustrated in Figure~\ref{fig:transitions} depending on whether $\tau_{i_{j+1}}$ is counterclockwise or clockwise from $\tau_{i_j}$. That is, if $\tau_{i_{j+1}}$ is counterclockwise to~$\tau_{i_j}$, the next steps of $\rho_{\gamma}$ will consist of a step two move that crosses $\tau_{i_j}$ and then a step one move from $\tau_{i_j}$ to $\tau_{i_{j+1}}$. If $\tau_{i_{j+1}}$ is clockwise from $\tau_{i_j}$, we cross $\tau_{i_j}$ with a step two move, then we reach $\tau_{i_{j+1}}$ through a step one move to reach $a_j$, then a step three move to travel along $a_j$ and we finish with a step one move to arrive at $\tau_{i_j}$.

After $(d-1)$ transitions, we reach the final triangle. We cross $\tau_{i_d}$ with a step three move, we then apply a step one move to get to $z$ and then we travel along $z$ with a step three move to reach the point $v^{\pm}_{Q, z}$. We call this the standard $M$-path associated to $\gamma$.

\textbf{Closed loops:} In the case where $\gamma$ is a closed loop, the process is largely the same. To begin, we choose a triangle $\Delta$ in $T$ such that two of its arcs are crossed by $\gamma$ and we label the two arcs $\tau_{i_1}$ and $\tau_{i_d}$ so that $\tau_{i_1}$ is clockwise from $\tau_{i_d}$. We label the third side of $\Delta$ by $a$. Let $p$ be a point on $\gamma$ which lies in $\Delta$ and has the form $v^{\pm}_{*, \tau_{i_1}}$. We then let $\tau_{i_1}, \ldots,\tau_{i_d}$ denote the ordered sequence of arcs which are crossed by $\gamma$, when one travels from $p$ away from $\Delta$. The standard $M$-path $\rho_{\gamma}$ associated to $\gamma$ is defined just like it was for generalized arcs, that is, by following the sequences of steps illustrated by Figure~\ref{fig:transitions}, depending on whether $\tau_{i_{j+1}}$ is counterclockwise or clockwise from $\tau_{i_j}$. The standard $M$-path begins and ends at the point $p$ and we consider the indices modulo $n$.
\end{Definition}

If $\rho_{\gamma}$ is not a standard $M$-path, for the generalized arc or loop $\gamma$, then one can deform $\gamma$ into a standard $M$-path $\rho_{\gamma}$ by making local adjustments that do not influence the upper right entry or trace.
We are now ready to prove Theorem~\ref{thm:matrixformulalamination}, which states that, with our additional adjustments, the matrix formulae of~\cite{MW13} hold for arbitrary principal laminations.

\begin{proof}[Proof of Theorem~\ref{thm:matrixformulalamination}]
 We first consider the case where $\gamma$ is a generalized arc and consider its standard $M$-path, $\rho_{\gamma}$. The first two steps of this path correspond to the matrix product
 \[\begin{bmatrix}
 1 & 0\\
 \dfrac{x_b}{x_a x_{i_1}} & 1
 \end{bmatrix}
 \begin{bmatrix}
 0 & x_a\\
 -\dfrac{1}{x_a} & 0
 \end{bmatrix}=
 \begin{bmatrix}
 0 & x_a\\
 -\dfrac{1}{x_a} & \dfrac{x_b}{x_{i_1}}
 \end{bmatrix},\]
 while the last three steps correspond to
\begin{align*}
\begin{bmatrix}
 0 & x_z\\
 -\dfrac{1}{x_z} & 0
 \end{bmatrix}
 \begin{bmatrix}
 1 & 0\\
 \dfrac{x_w}{x_{i_d} x_z} & 1
 \end{bmatrix}
 \begin{bmatrix}
 1 & 0\\
 0 & y_{i_d}
 \end{bmatrix}&=
 \begin{bmatrix}
 \dfrac{x_w}{x_{i_d}} & x_z y_{i_d}\vspace{1mm}\\
 -\dfrac{1}{x_z} & 0
 \end{bmatrix}, \qquad\text{when}\quad b_T(\mathcal{L}_T, i_d)=1,\\
\begin{bmatrix}
 0 & x_z\\
 -\frac{1}{x_z} & 0
 \end{bmatrix}
 \begin{bmatrix}
 1 & 0\\
 \dfrac{x_w}{x_{i_d} x_z} & 1
 \end{bmatrix}
 \begin{bmatrix}
 y_{i_d} & 0\\
 0 & 1
 \end{bmatrix}&=
 \begin{bmatrix}
 \dfrac{x_w y_{i_d}}{x_{i_d}} & x_z\vspace{1mm}\\
 -\dfrac{y_{i_d}}{x_z} & 0
 \end{bmatrix},\qquad \text{when}\quad b_T(\mathcal{L}_T, i_d)=-1.
\end{align*}

In between the first and the last steps, for the portion between $\tau_{i_j}$ and $\tau_{i_{j+1}}$, where $1\leq j\leq d-1$, we have the following four cases that depend on the sign of $b_T(\mathcal{L}_T, \tau_{i_j})$ and whether $\tau_{i_{j+1}}$ lies clockwise or counterclockwise from $\tau_{i_j}$,{\samepage
\begin{gather*}
\begin{split}
 &\begin{bmatrix}
 1 & 0\\
 \dfrac{x_{a_j}}{x_{i_j}x_{i_{j+1}}} & 1
 \end{bmatrix}
 \begin{bmatrix}
 1 & 0\\
 0 & y_{i_j}
 \end{bmatrix}=
 \begin{bmatrix}
 1 & 0\\
 \dfrac{x_{a_j}}{x_{i_j}x_{i_{j+1}}} & y_{i_j}
 \end{bmatrix} , \qquad\text{when}\quad b_T(\mathcal{L}_T, i_j)=1, \text{CCW},\\
& \begin{bmatrix}
 1 & 0\\
 \dfrac{x_{a_j}}{x_{i_j}x_{i_{j+1}}} & 1
 \end{bmatrix}
 \begin{bmatrix}
 y_{i_j} & 0\\
 0 & 1
 \end{bmatrix}=
 \begin{bmatrix}
 y_{i_j} & 0\vspace{1mm}\\
 \dfrac{x_{a_j}y_{i_j}}{x_{i_j}x_{i_{j+1}}} & 1
 \end{bmatrix}, \qquad\text{when}\quad b_T(\mathcal{L}_T, i_j)=-1, \text{CCW},\\
& \begin{bmatrix}
 1 & 0\\
 \dfrac{x_{i_j}}{x_{a_j}x_{i_{j+1}}} & 1\\
 \end{bmatrix}
 \begin{bmatrix}
 0 & x_{a_j}\\
 -\dfrac{1}{x_{a_j}} & 0
 \end{bmatrix}
 \begin{bmatrix}
 1 & 0\\
 \dfrac{x_{i_{j+1}}}{x_{a_j}x_{i_j}} & 1
 \end{bmatrix}
 \begin{bmatrix}
 1 & 0\\
 0 & y_{i_j}
 \end{bmatrix}=
 \begin{bmatrix}
 \dfrac{x_{i_{j+1}}}{x_{i_j}} & x_{a_j}y_{i_j}\\
 0 & \dfrac{x_{i_j}y_{i_j}}{x_{i_{j+1}}}
 \end{bmatrix}, \\
 &\qquad \text{when}\quad b_T(\mathcal{L}_T, i_j)=1, \text{CW},\\
& \begin{bmatrix}
 1 & 0\\
 \dfrac{x_{i_j}}{x_{a_j}x_{i_{j+1}}} & 1\\
 \end{bmatrix}
 \begin{bmatrix}
 0 & x_{a_j}\\
 -\dfrac{1}{x_{a_j}} & 0
 \end{bmatrix}
 \begin{bmatrix}
 1 & 0\\
 \dfrac{x_{i_{j+1}}}{x_{a_j}x_{i_j}} & 1
 \end{bmatrix}
 \begin{bmatrix}
 y_{i_j} & 0\\
 0 & 1
 \end{bmatrix}=
 \begin{bmatrix}
 \dfrac{x_{i_{j+1}}y_{i_j}}{x_{i_j}} & x_{a_j}\\
 0 & \dfrac{x_{i_j}}{x_{i_{j+1}}}
 \end{bmatrix} \\
& \qquad\text{when}\quad b_T(\mathcal{L}_T, i_j)=-1, \text{CW}.
\end{split}
\end{gather*}
Here, $a_j$ is the third side of the triangle with sides $\tau_{i_j}$, $\tau_{i_{j+1}}$.}

By Definition~\ref{definition:chilamination} and Theorem~\ref{theorem:perfectmatchingenum}, we find
 \begin{itemize}\itemsep=0pt
 \item {\samepage\[
 \chi_{\gamma, \mathcal{L}_{T}}=\operatorname{ur}\left(\begin{bmatrix}
\dfrac{x_w}{x_{i_d}} && x_{z}y_{i_d}\vspace{1mm}\\
-\dfrac{1}{x_{z}} && 0
\end{bmatrix}
M
\begin{bmatrix}
0 && x_a\\
-\dfrac{1}{x_a} && \dfrac{x_b}{x_{i_1}}
\end{bmatrix}\right)=\frac{\sum_P x(P)y_{\mathcal{L}_T}(P)}{\operatorname{cross}_T(\gamma)},\]
when $ b_{T}(\mathcal{L}_{T}, i_{d})=1$,}

\item \[
\chi_{\gamma, \mathcal{L}_{T}}=\operatorname{ur}\left(\begin{bmatrix}
\dfrac{x_w y_{i_d}}{x_{i_d}} && x_{z}\vspace{1mm}\\
-\dfrac{y_{i_d}}{x_{z}} && 0
\end{bmatrix}
M
\begin{bmatrix}
0 && x_a\\
-\dfrac{1}{x_a} && \dfrac{x_b}{x_{i_1}}
\end{bmatrix}\right)=\frac{\sum_P x(P)y_{\mathcal{L}_T}(P)}{\operatorname{cross}_T(\gamma)},\]
 when $ b_{T}(\mathcal{L}_{T}, i_{d})=-1$,
\end{itemize}
where the middle matrix $M$ is obtained by multiplying together a sequence of matrices of the four forms mentioned previously. The term $\chi_{\gamma, \mathcal{L}_T}$ has an interpretation in terms of perfect matchings of some abstract snake graph $G$. Namely, the abstract snake graph $G$ is precisely~$\mathcal{G}_{T, \gamma}$ associated to $\gamma$, meaning $M=M_d$ from Definition~\ref{definition:tilematrix}, which justifies the right-most equality and the proof is complete.

The case where $\gamma$ is a closed loop is virtually the same, the only difference being that our initial and final steps are different. By construction of the standard $M$-path, we have that $\tau_{i_1}$ is in the clockwise direction of $\tau_{i_d}$. This implies the final steps of $\rho_{\gamma}$ corresponds to either
\begin{align*}
 &\begin{bmatrix}
 \dfrac{x_{i_1}}{x_{i_d}} && x_a y_{i_d}\\
 0 && \dfrac{x_{i_d}y_{i_d}}{x_{i_1}}
 \end{bmatrix} ,\qquad\text{when}\quad b_T(\mathcal{L}_T, \gamma)=1,\\
 &\begin{bmatrix}
 \dfrac{x_{i_1}y_{i_d}}{x_{i_d}} && x_a\\
 0 && \dfrac{x_{i_d}}{x_{i_1}}
 \end{bmatrix} ,\qquad\text{when}\quad b_T(\mathcal{L}_T, \gamma)=-1.
\end{align*}

By Definition~\ref{definition:chilamination}, Theorem~\ref{theorem:goodmatchingenum} and using the interpretation of $\chi_{\gamma, \mathcal{L}_T}$ in terms of good matchings of some abstract band graph $G$, we obtain
\begin{itemize}\itemsep=0pt
 \item \[
 \chi_{\gamma, \mathcal{L}_{T}}=\operatorname{tr}\left(\begin{bmatrix}
\dfrac{x_{i_1}}{x_{i_d}} & x_{a}y_{i_d}\\
0 & \dfrac{x_{i_d}y_{i_d}}{x_{i_1}}
\end{bmatrix}
M\right)=\frac{\sum_P x(P)y_{\mathcal{L}_T}(P)}{\operatorname{cross}_T(\gamma)},
\]
 when $ b_{T}(\mathcal{L}_{T}, i_{d})=1$,

\item \[
\chi_{\gamma, \mathcal{L}_{T}}=\operatorname{tr}\left(\begin{bmatrix}
\dfrac{x_{i_1}y_{i_d}}{x_{i_d}} & x_{a}\\
0 & \dfrac{x_{i_d}}{x_{i_1}}
\end{bmatrix}
M\right)=\frac{\sum_P x(P)y_{\mathcal{L}_T}(P)}{\operatorname{cross}_T(\gamma)},\]
 when $ b_{T}(\mathcal{L}_{T}, i_{d})=-1$,
\end{itemize}
where again $M$ is a obtained by multiplying together a sequence of matrices of the four forms mentioned previously and can be interpreted as the matrix $M_d$ for $\overline{\mathcal{G}}_{T, \gamma}$.
\end{proof}

\subsection{Example with closed curve}
\begin{figure}[!ht]
 \centering
 \includegraphics[scale=.55]{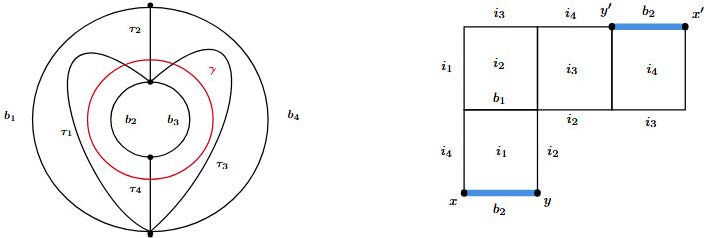}
 \caption{Triangulation $T$ and band graph $\widetilde{\mathcal{G}}_{T, \gamma}$ for the closed loop $\gamma$.}\label{fig:MWBand}
\end{figure}

\begin{figure}[!ht]
 \centering
 \includegraphics[scale=.48]{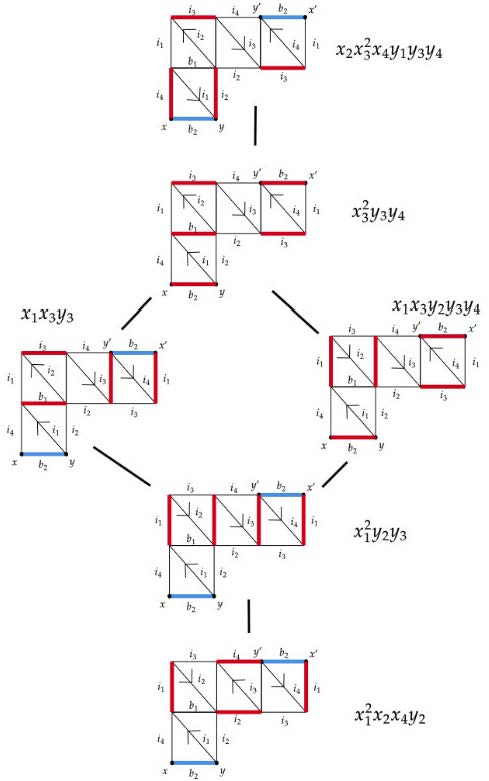}
 \caption{Perfect matchings on $\widetilde{\mathcal{G}}_{T, \gamma}$, along with their associated weight monomial. Note, $y_{\mathcal{L}_T'}(P)$ is calculated using Definition~\ref{definition:L_Toriented} for each perfect matching~$P$.}\label{fig:bandperfectmatchings}
\end{figure}

Consider the triangulation and closed loop $\gamma$ with band graph $\widetilde{\mathcal{G}}_{T, \gamma}$ pictured in Figure~\ref{fig:MWBand}.
We~compute $\chi_{\gamma, \mathcal{L}_T'}$ as well as
\[
\frac{\sum_P x(P)y_{\mathcal{L}_T'}(P)}{\operatorname{cross}_T(\gamma)}\]
 and show they are equal. The perfect matchings and how they relate to one another are pictured in Figure~\ref{fig:bandperfectmatchings}. Taking the sum of each monomial in the Figure~\ref{fig:bandperfectmatchings}, we find \[\frac{\sum_P x(P)y_{\mathcal{L}_T'}(P)}{\operatorname{cross}_T(\gamma)}=\frac{y_2x_1^2x_2x_4+y_2y_3x_1^2+(y_3+y_2y_3y_4)x_1x_3+y_3y_4x_3^2+y_1y_3y_4x_2x_3^2x_4}{x_1x_2x_3x_4}.\]

Meanwhile, the matrix formula yields
\begin{align*}
\chi_{\gamma, \mathcal{L}_T'}={}&\operatorname{tr}\left(
\begin{bmatrix}
 \dfrac{x_1}{x_4} & y_4b_2\\
 0 & \dfrac{x_4y_4}{x_1}
\end{bmatrix}
\begin{bmatrix}
 \dfrac{x_4}{x_3} & y_3b_3\\
 0 & \dfrac{x_3y_3}{x_4}
\end{bmatrix}
\begin{bmatrix}
 y_2 & 0\\
 \dfrac{b_4y_2}{x_2x_3} & 1
\end{bmatrix}
\begin{bmatrix}
 1 & 0\\
 \dfrac{b_1}{x_1x_2} & y_1
\end{bmatrix}\right)\\
={}&\operatorname{tr}
\begin{bmatrix}
 \dfrac{\begin{array}{@{}c@{}}y_2y_3y_4x_1x_3b_2b_4+y_3y_4x_3^2b_1b_2+y_2y_3x_1^2b_3b_4\\
 +y_2x_1^2x_2x_4+y_3x_1x_3b_1b_3\end{array}}{x_1x_2x_3x_4} & \dfrac{y_1y_3y_4x_3b_2+y_1y_3x_1b_3}{x_4}\vspace{2mm}\\
 \dfrac{y_2y_3y_4x_1b_4+y_3y_4x_3b_1}{x_1^2x_2} & \dfrac{y_1y_3y_4x_3}{x_1}
\end{bmatrix}\\
={}&\dfrac{y_1y_3y_4x_2x_3^2x_4+y_2y_3y_4x_1x_3b_2b_4+y_3y_4x_3^2b_1b_2+y_2y_3x_1^2b_3b_4+y_2x_1^2x_2x_4}{x_1x_2x_3x_4}\\
&+\dfrac{y_3x_1x_3b_1b_3}{x_1x_2x_3x_4}.
\end{align*}
Letting $b_1=b_2=b_3=b_4=1$ and rearranging the terms in our numerator, we conclude
\begin{align*}
\chi_{\gamma, \mathcal{L}_T'}&=\dfrac{\sum_P x(P)y_{\mathcal{L}_T'}(P)}{\operatorname{cross}_T(\gamma)}\\
&=\dfrac{y_2x_1^2x_2x_4+y_2y_3x_1^2+(y_3+y_2y_3y_4)x_1x_3+y_3y_4x_3^2+y_1y_3y_4x_2x_3^2x_4}{x_1x_2x_3x_4}.\end{align*}

\begin{Remark}[{comparison to \cite[Example~3.23]{MW13}}] We chose to revisit \cite[Example 3.23]{MW13} using a~different principal lamination to demonstrate the differences in our construction. In \cite{MW13}, they use the same triangulation and closed curve with the principal lamination $\mathcal{L}_T$ defined by~$b_T(\mathcal{L}_T, \tau_{i})=1$ for each $1\leq i\leq 4$. For our example, we consider $\mathcal{L}_T'$ defined by taking $\mathcal{L}_T$ and changing~$b_T(\mathcal{L}_T, \tau_2)$ from $1$ to $-1$ and keeping everything else the same.

With these differences, Musiker and Williams obtain
\begin{align*}
\chi_{\gamma, \mathcal{L}_T}&=\frac{\sum_P x(P)y_{\mathcal{L}_T}(P)}{\operatorname{cross}_T(\gamma)}\\
&=\frac{x_1^2x_2x_4+y_3x_1^2+(y_2y_3+y_3y_4)x_1x_3+y_2y_3y_4x_3^2+y_1y_2y_3y_4x_2x_3^2x_4}{x_1x_2x_3x_4},\end{align*}
whereas, we have
\begin{align*}
\chi_{\gamma, \mathcal{L'}_T}&=\frac{\sum_P x(P)y_{\mathcal{L'}_T}(P)}{\operatorname{cross}_T(\gamma)}\\
&=\frac{y_2x_1^2x_2x_4+y_2y_3x_1^2+(y_3+y_2y_3y_4)x_1x_3+y_3y_4x_3^2+y_1y_3y_4x_2x_3^2x_4}{x_1x_2x_3x_4}.\end{align*}
\end{Remark}

\begin{Remark}
Also, notice that Figure~\ref{fig:bandperfectmatchings} is arranged to look like the Hasse diagram for a~poset. Although we may always arrange the set of (good) perfect matchings as a poset with the cover relation given by ``flipping'' local tiles as in \cite{MSW11}, the coefficient variables cannot be directly recovered without the lamination. For instance, the minimal element in this poset has a non-trivial height monomial of $y_2$. This is not typically how the $y$-coefficients or height monomial is classically defined. In general, there will exist a choice of isotopic representative of arc that gives the classical poset structure interpretation of the coefficients, but perturbations of this arc will change the poset structure. An example of this can be seen comparing the poset we obtain in Figures~\ref{fig:posetforalpha'} and~\ref{fig:posetgraph} in Section~\ref{section:expansionformulanonorientable}.
\end{Remark}

\section[Expansion formula for one-sided closed curves using matrix products]{Expansion formula for one-sided closed curves\\ using matrix products}\label{section:expansionformulanonorientable}

In this section, we will prove the matrix formulae in the previous section can be used to find the Laurent expansion for one-sided closed curves. To prove this result, we will be following the general guidelines found in Section~\ref{section:matrixorientable} of this paper. That is, we will be associating $2\times 2$ matrices to the parallelograms of a band graph, which will then give us a way of representing the good matching enumerator of our graph in terms of the trace of a product of matrices. To conclude this section, we will show that the matrix product associated with the band graphs of one-sided curves coincides with the matrix product associated to a canonical $M$-path associated to the one-sided curve.

\subsection{Good matching enumerators for one-sided closed curves}

The basic ingredients used to construct band graphs for one-sided curves are largely the same as those found in Definition~\ref{definition:bandgraph}, but there is a slight twist.

For two-sided closed curves, we always glue together sides with the same sign, but when working with the one-sided closed curve, we will now be identifying edges that have the opposite sign. With this in mind, the following abstract band graph will differ from the one defined earlier.

\begin{Definition}[abstract band graph for one-sided curves]\label{definition:quasibandgraph}
An abstract band graph for a~one-sided curve with $d$ tiles is formed by concatenating the following puzzle pieces:
\begin{itemize}\itemsep=0pt
\item an initial triangle

\begin{figure}[!ht]
 \centering
 \includegraphics[scale=.5]{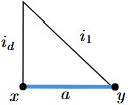}
\end{figure}

\item $d-1$ parallelograms $H_1, \ldots, H_{d-1}$, where each $H_j$ is as before,

\item a final triangle based on whether $d$ is odd or even.

\begin{figure}[!ht]
 \centering
 \includegraphics[scale=.5]{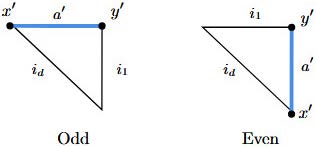}
\end{figure}

\end{itemize}
\end{Definition}

Just like the orientable case, we associate various matrices to the parallelograms of our graph that will be determined by comparing the shape of the parallelograms $H_i$ and $H_{i-1}$, as well as the principal lamination used for coefficients. See Definition~\ref{definition:tilematrix}.

We now provide a description of the entries of $M_d$ in terms of weight and coefficient monomials~$x(P)$ and $y_{\mathcal{L}_T}(P)$ for good matchings of our new variant of abstract band graphs. The following proposition is an immediate corollary of Proposition~\ref{prop:snakematrix} after relabeling the edges of our graph.

\begin{Proposition}\label{prop:quasibandmatrix}
Let $G$ be an abstract snake graph with $d$-tiles, but with a labeling obtained by substituting $w$ with $a'$, $b$ with $i_d$ and $z$ with $i_1$. Write \smash{$M_d=\big[\begin{smallmatrix}
A_d & B_d\\
C_d & D_d
\end{smallmatrix}\big]$}.
For $d\geq 2$,
\begin{itemize}\itemsep=0pt
\item when $b_{T}(\mathcal{L}_{T}, {i_{d-1}})=1$, we have
\begin{alignat*}{3}
& A_d = \frac{\sum_{P\in S_A}x(P)y_{\mathcal{L}_T}(P)}{(x_{i_1}\cdots x_{i_{d-1}})x_ax_{a'}} ,\qquad&& B_d = \frac{\sum_{P\in S_B}x(P)y_{\mathcal{L}_T}(P)}{(x_{i_2}\cdots x_{i_{d-1}})x_{i_d}x_{a'}},\\
&C_d = \frac{\sum_{P\in S_C}x(P)y_{\mathcal{L}_T}(P)}{(x_{i}\cdots x_{i_d})x_ax_{i_1}y_{i_d}} ,\qquad&& D_d = \frac{\sum_{P\in S_D}x(P)y_{\mathcal{L}_T}(P)}{(x_{i_2}\cdots x_{i_{d}})x_{i_d}x_{i_1}y_{i_d}},&
\end{alignat*}
\item when $b_{T}(\mathcal{L}_{T}, i_{d-1})=-1$, we have
\begin{alignat*}{3}
& A_d = \frac{\sum_{P\in S_A}x(P)y_{\mathcal{L}_T}(P)}{(x_{i_1}\cdots x_{i_{d-1}})x_ax_a'y_{i_d}} ,\qquad&& B_d = \frac{\sum_{P\in S_B}x(P)y_{\mathcal{L}_T}(P)}{(x_{i_2}\cdots x_{i_{d-1}})x_{i_d}x_{a'} y_{i_d}},&\\
& C_d = \frac{\sum_{P\in S_C}x(P)y_{\mathcal{L}_T}(P)}{(x_{i}\cdots x_{i_d})x_ax_{i_1}},\qquad&& D_d = \frac{\sum_{P\in S_D}x(P)y_{\mathcal{L}_T}(P)}{(x_{i_2}\cdots x_{i_{d}})x_{i_d}x_{i_1}}.&
\end{alignat*}
\end{itemize}
\noindent For $d=1$, the formulae remain the same, but the cases are determined by the sign of $b_T(\mathcal{L}_T, i_1)$.

\noindent Here, $S_A$, $S_B$, $S_C$ and $S_D$ are the sets of perfect matchings of $G$ which use the edges $\{a, a'\}$, $ \{i_d, a'\}$, $\{a, i_1\}$ and $\{i_d, i_1\}$, respectively.
\end{Proposition}

The following comes as a direct corollary of Proposition~\ref{prop:quasibandmatrix}.

\begin{Corollary}\label{cor:matrixformula}
Let $G$ be an abstract band graph corresponding to a one-sided closed curve with~$d$ tiles. Write \smash{$M_d=\big[\begin{smallmatrix}
 A_d & B_d\\
C_d & D_d
\end{smallmatrix}\big]$}.
\begin{itemize}\itemsep=0pt
\item When $b_T({\mathcal{L}_T}, i_d)=1$, then
\begin{align}
\frac{\sum_P x(P)y(P)}{x_{i_1}\cdots x_{i_d}}&=\frac{x_a A_d}{x_{i_d}}+\frac{B_d}{x_{i_1}}+x_{i_1}y_{i_d}C_d\nonumber\\
&=\operatorname{tr}\left(\begin{bmatrix}
\dfrac{x_a}{x_{i_d}} & x_{i_1}y_{i_d}\vspace{1mm}\\
\dfrac{1}{x_{i_1}} & 0
\end{bmatrix}
\begin{bmatrix}
A_d & B_d\\
C_d & D_d
\end{bmatrix}\right).\label{eqn:matrixformula1+}
\end{align}

\item When $b_T({\mathcal{L}_T}, i_d)=-1$, then
\begin{align}
\frac{\sum_P x(P)y(P)}{x_{i_1}\cdots x_{i_d}}&=\frac{x_ay_{i_d} A_d}{x_{i_d}}+\frac{y_{i_d}B_d}{x_{i_1}}+x_{i_1}C_d\nonumber\\
&=\operatorname{tr}\left(\begin{bmatrix}
\dfrac{x_ay_{i_d}}{x_{i_d}} && x_{i_1}\\
\dfrac{y_{i_d}}{x_{i_1}} && 0
\end{bmatrix}
\begin{bmatrix}
A_d && B_d\\
C_d && D_d
\end{bmatrix}\right), \label{eqn:matrixformula2+}
\end{align}
where the sum is over all good matchings of $G$.
\end{itemize}
\end{Corollary}
\begin{proof}
Following the notation from Proposition~\ref{prop:quasibandmatrix}, we consider the sets $S_A$, $S_B$, $S_C$ and $S_D$. If $G$ is our snake graph from Proposition~\ref{prop:quasibandmatrix} and $\widetilde{G}$ is the band obtained from identifying $a$ and~$a'$, then every perfect matching from $S_A$, $S_B$ and $S_C$ descends to a good matching of $\widetilde{G}$ after removing either $a$ or $a'$. On the other hand, no perfect matching from $S_D$ descends to a~good matching of $\widetilde{G}$. As all good matchings of $\widetilde{G}$ are obtained uniquely from perfect matchings in $S_A$, $S_B$ or $S_C$, the result follows.
\end{proof}

The corollary above immediately proves the following theorem.

\begin{Theorem}\label{theorem:quasigoodmatchingenum}
Suppose $G$ is an abstract band graph with $d$ tiles that corresponds to a one-sided closed curve. Then its good matching enumerator is given by
\begin{itemize}\itemsep=0pt
\item \[
\sum_Px(P)y(P)=x_{i_1}\ldots x_{i_d} \operatorname{tr}\left(\begin{bmatrix}
\frac{x_a}{x_{i_d}} & x_{i_1}y_{i_d}\\
\frac{1}{x_{i_1}} & 0
\end{bmatrix}M_d\right),\]
when $ b_{T}(\mathcal{L}_{T}, i_{d})=1$,
\item \[\sum_Px(P)y(P)=x_{i_1}\ldots x_{i_d} \operatorname{tr}\left(\begin{bmatrix}
\frac{x_ay_{i_d}}{x_{i_d}} & x_{i_1}\\
\frac{y_{i_d}}{x_{i_1}} & 0
\end{bmatrix}M_d\right),\]
 when $ b_{T}(\mathcal{L}_{T}, i_{d})=-1$,
\end{itemize}
where the sum is over all good matchings of $G$.
\end{Theorem}

\subsection[The standard M-path for a one-sided closed curve]{The standard $\boldsymbol{M}$-path for a one-sided closed curve}

In this subsection, for any one-sided closed curve $\alpha$, we will associate to it a standard $M$\nobreakdash-path~$\rho_{\alpha}$ and show that the associated matrix formula will have the same form as equations \eqref{eqn:matrixformula1+} and~\eqref{eqn:matrixformula2+} found in Corollary~\ref{cor:matrixformula}.

The elementary steps used in the standard $M$-path will be
the same as those used in the orientable case; however,
we will need to make an adjustment to the matrix corresponding to the type III step that
travels through the crosscap, as we will now be
traveling along an arc from a point of the form
$v_{m,\tau}^{+}$ ($v_{m,\tau}^{-}$) to another point of
the form $v_{m', \tau}^{+}$ ($v_{m',\tau}^{-}$).

\begin{Definition}[elementary step of type 3$'$]\label{definition:newstep}
For this variation of the third type of elementary step, we travel through the crosscap along a path parallel to a fixed arc $\tau$ connecting two points~\smash{$v^{\pm}_{m, \tau}$} and \smash{$v^{\pm}_{m', \tau}$} associated to distinct marked points $m$, $m'$. The associated matrix is~\smash{$\big[\begin{smallmatrix}
0 & \pm x_{\tau}\\
\pm \frac{1}{x_{\tau}} & 0
\end{smallmatrix}\big]$}, where we use $x_r$, $\frac{1}{x_r}$ if this step sees $\tau$ on the right and $-x_r$, $\frac{-1}{x_r}$ otherwise.
\end{Definition}

\begin{figure}[!ht]
\centering
\includegraphics[scale=.49]{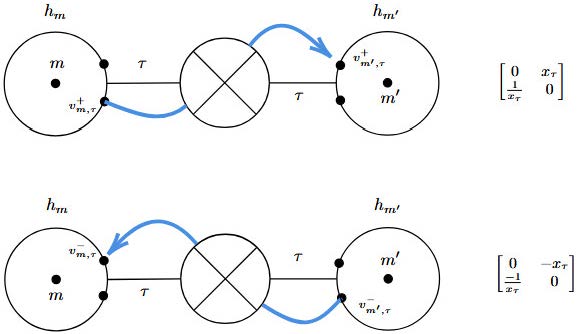}

\caption{Elementary steps of type 3$'$ in the positive and negative direction.}
\end{figure}

An illustration of the type 3$'$ elementary step can be found in Figure~\ref{fig:mobiuspath} of Section~\ref{section:example}. The next lemma is a quick check verification that Lemma~\ref{lemma: independenceofpath} still holds with our new type three step.

\begin{Lemma}\label{lemma:independenceofpath2}
Fix $\SM$ and $T$. Let $\alpha$ be a one-sided closed curve. Then for any two $M$-paths, $\rho$ and $\rho'$ associated to $\alpha$, we have $|{\operatorname{tr}(M(\rho))}|=|{\operatorname{tr}(M(\rho'))}|$.
\end{Lemma}
\begin{proof}
First, we observe if we have an $M$-path $M(\rho)=M(\rho_n)\cdots M(\rho_1)$, then \[\operatorname{tr}[M(\rho_n)\cdots M(\rho_1)]=\operatorname{tr}[M(\rho_{i})\cdots M(\rho_{1})\cdot M(\rho_n)\cdots M(\rho_{i+1})]\] for all $1\leq i\leq n$, so the trace is invariant under the starting point of our $M$-path. The original three elementary steps are covered already by Lemma \ref{lemma: independenceofpath}.

Furthermore, we have the equality \[
\begin{bmatrix}
 0 & x_{\tau}\\
 \dfrac{1}{x_{\tau}} & 0
\end{bmatrix}
\begin{bmatrix}
 y_{\tau} & 0\\
 0 & 1
\end{bmatrix}=\begin{bmatrix}
 1 & 0\\
 0 & y_{\tau}
\end{bmatrix}\begin{bmatrix}
 0 & x_{\tau}\\
 \dfrac{1}{x_{\tau}} & 0
\end{bmatrix},\]
which implies, with respect to $|{\operatorname{tr}(M(\rho))}|$, crossing $\tau$ with a step two move and then passing through the crosscap with the adjusted step three move is the same as first passing through the crosscap and then crossing $\tau$.
\end{proof}

Despite working on a non-orientable surface, the algorithm for constructing the standard $M$-path for the one-sided closed curve will be reminiscent to the algorithm for creating standard $M$-paths for closed curves on orientable surfaces. The differences being that we require the first arc to be counter-clockwise from the last; and we require that the final step of the path corresponds to traveling through the crosscap.

\begin{Definition}[standard $M$-path for a one-sided closed curve]\label{definition:Mpath2}
Let $\alpha$ be a one-sided closed curve which crosses $d$ arcs of a fixed triangulation $T$ (counted with multiplicity). Choose a~triangle~$\Delta$ in $T$ such that two of its arcs are crossed by $\alpha$. We label these two arcs as $\tau_{i_1}$ and~$\tau_{i_d}$, where $\tau_{i_1}$ is in the counterclockwise direction from $\tau_{i_d}$. We label the third side with $a$. Let $p$ be a point on $\alpha$ which lies in $\Delta$ and has the form $v_{*, \tau_{i_1}}^{\pm}$. Letting, $\tau_{i_1},\ldots, \tau_{i_d}$ denote the ordered sequence of arcs that are crossed by $\alpha$ as we move from $p$ and $\Delta$, we have that the standard $M$-path $\rho_{\alpha}$ associated to $\alpha$ is the same as Definition~\ref{definition:Mpath}, starting and ending at $p$ and traveling along elementary steps based on whether $\tau_{i_{j+1}}$ is counterclockwise or clockwise from $\tau_{i_j}$. The final elementary steps of $\rho_{\alpha}$ will be determined by the fact that $\tau_{i_1}$ is in the counterclockwise direction from $\tau_{i_d}$ and that, to return to $p$, we must travel parallel to $i_1$ through the crosscap. We consider indices modulo $n$. See Figure $\ref{fig:mobiuspath}$ for an example of a such a standard $M$-path.
\end{Definition}

We are ready to prove the main theorem of this section.

\begin{Theorem}\label{thm:matrixformula}
Let $\SM$ be a marked non-orientable surface with a triangulation $T$, and let $T=\{\tau_{i_1}, \ldots, \tau_{i_n}\}$ be the corresponding triangulation. Let $\mathcal{A}_{\mathcal{L}_T}\SM$ be the corresponding quasi-cluster algebra associated to the principal lamination $\mathcal{L}_T$.

Suppose $\alpha$ is a one-sided closed curve which is not contractible, and has no contractible kinks. Then{\samepage \[\displaystyle\chi_{\alpha, \mathcal{L_T}}=\frac{1}{\operatorname{cross}_T(\alpha)}\sum_{P}x(P)y_{L_T}(P),\] where the sum is over all good matchings $P$ of the band graph $\widetilde{\mathcal{G}}_{T, \alpha}$.}
\end{Theorem}
\begin{proof}
We consider the standard $M$-path of $\alpha$ defined in Definition~\ref{definition:Mpath2}. For the indices $j$ from~$\tau_{1}$ to $\tau_{d-1}$, we have the following four cases that depend on the sign of $b_T(\mathcal{L}_T, \tau_{i_j})$ and whether~$\tau_{i_{j+1}}$ lies clockwise or counterclockwise from $\tau_{i_j}$:
\begin{gather}
 \begin{bmatrix}
 1 & 0\\
 \dfrac{x_{a_j}}{x_{i_j}x_{i_{j+1}}} & 1
 \end{bmatrix}
 \begin{bmatrix}
 1 & 0\\
 0 & y_{i_j}
 \end{bmatrix}=
 \begin{bmatrix}
 1 & 0\\
 \dfrac{x_{a_j}}{x_{i_j}x_{i_{j+1}}} & y_{i_j}
 \end{bmatrix},\qquad\text{when}\quad b_T(\mathcal{L}_T, i_j)=1, \text{CCW},\label{eqn:ccw+1}\\
 \begin{bmatrix}
 1 & 0\\
 \dfrac{x_{a_j}}{x_{i_j}x_{i_{j+1}}} & 1
 \end{bmatrix}
 \begin{bmatrix}
 y_{i_j} & 0\\
 0 & 1
 \end{bmatrix}=
 \begin{bmatrix}
 y_{i_j} & 0\vspace{1mm}\\
 \dfrac{x_{a_j}y_{i_j}}{x_{i_j}x_{i_{j+1}}} & 1
 \end{bmatrix} ,\qquad\text{when}\quad b_T(\mathcal{L}_T, i_j)=-1, \text{CCW},\label{eqn:ccw-1}\\
 \begin{bmatrix}
 1 & 0\\
 \dfrac{x_{i_j}}{x_{a_j}x_{i_{j+1}}} & 1\\
 \end{bmatrix}
 \begin{bmatrix}
 0 & x_{a_j}\\
 -\dfrac{1}{x_{a_j}} & 0
 \end{bmatrix}
 \begin{bmatrix}
 1 & 0\\
 \dfrac{x_{i_{j+1}}}{x_{a_j}x_{i_j}} & 1
 \end{bmatrix}
 \begin{bmatrix}
 1 & 0\\
 0 & y_{i_j}
 \end{bmatrix}=
 \begin{bmatrix}
 \dfrac{x_{i_{j+1}}}{x_{i_j}} & x_{a_j}y_{i_j}\\
 0 & \dfrac{x_{i_j}y_{i_j}}{x_{i_{j+1}}}
 \end{bmatrix} ,\nonumber\\
 \qquad \text{when}\quad b_T(\mathcal{L}_T, i_j)=1, \text{CW},\label{eqn:cw+1}\\
 \begin{bmatrix}
 1 & 0\\
 \dfrac{x_{i_j}}{x_{a_j}x_{i_{j+1}}} & 1\\
 \end{bmatrix}
 \begin{bmatrix}
 0 & x_{a_j}\\
 -\dfrac{1}{x_{a_j}} & 0
 \end{bmatrix}
 \begin{bmatrix}
 1 & 0\\
 \dfrac{x_{i_{j+1}}}{x_{a_j}x_{i_j}} & 1
 \end{bmatrix}
 \begin{bmatrix}
 y_{i_j} & 0\\
 0 & 1
 \end{bmatrix}=
 \begin{bmatrix}
 \dfrac{x_{i_{j+1}}y_{i_j}}{x_{i_j}} & x_{a_j}\\
 0 & \dfrac{x_{i_j}}{x_{i_{j+1}}}
 \end{bmatrix}, \nonumber\\
\qquad \text{when}\quad b_T(\mathcal{L}_T, i_j)=-1, \text{CW}.\label{eqn:cw-1}
\end{gather}
Here, $a_j$ is the third side of the triangle with sides $\tau_{i_j}$, $\tau_{i_{j+1}}$.

As $\tau_{i_1}$ is in the counterclockwise direction from $\tau_{i_{d}}$, and we must travel parallel along $\tau_{i_1}$ through the crosscap in order to return to $p$, the final few steps of $\rho_{\alpha}$ are represented by one of the two matrices
\begin{gather}
\begin{bmatrix}
0 & x_{i_1}\\
\dfrac{1}{x_{i_1}} & 0
\end{bmatrix}
\begin{bmatrix}
1 & 0\\
\dfrac{x_{a}}{x_{i_1}x_{i_d}} & y_{i_d}
\end{bmatrix}=\begin{bmatrix}
\dfrac{x_a}{x_{i_d}} & x_{i_1}y_{i_d}\vspace{1mm}\\
\dfrac{1}{x_{i_1}} & 0
\end{bmatrix}\qquad\text{when}\quad b_T(\mathcal{L}_T, i_d)=1, \label{eqn:finalstep1}\\
\begin{bmatrix}
0 & x_{i_1}\vspace{1mm}\\
\dfrac{1}{x_{i_1}} & 0
\end{bmatrix}
\begin{bmatrix}
y_{i_d} & 0\vspace{1mm}\\
\dfrac{x_{a}y_{i_d}}{x_{i_1}x_{i_d}} & 1
\end{bmatrix}=\begin{bmatrix}
\dfrac{x_ay_{i_d}}{x_{i_d}} & x_{i_1}\vspace{1mm}\\
\dfrac{y{_{i_d}}}{x_{i_1}} & 0
\end{bmatrix}\qquad\text{when}\quad b_T(\mathcal{L}_T, i_d)=-1. \label{eqn:finalstep2}
\end{gather}
By Definition~\ref{definition:chilamination}, Theorem~\ref{theorem:quasigoodmatchingenum} and using the interpretation of $\chi_{\alpha, \mathcal{L}_T}$ in terms of good matchings on some abstract band graph $\mathcal{G}$, we have
\begin{itemize}\itemsep=0pt
\item
\[
\chi_{\alpha, \mathcal{L}_T}=\operatorname{tr}\left(\begin{bmatrix}
\dfrac{x_a}{x_{i_d}} & x_{i_1}y_{i_d}\vspace{1mm}\\
\dfrac{1}{x_{i_1}} & 0
\end{bmatrix}M\right)=\frac{\sum_P x(P)y_{\mathcal{L}_T}(P)}{\operatorname{cross}_T(\alpha)},\]
when $ b_{T}(\mathcal{L}_{T}, i_{d})=1$,
\item \[
\chi_{\alpha, \mathcal{L}_T}=\operatorname{tr}\left(\begin{bmatrix}
\dfrac{x_ay_{i_d}}{x_{i_d}} & x_{i_1}\vspace{1mm}\\
\dfrac{y_{i_d}}{x_{i_1}} & 0
\end{bmatrix}M\right)=\frac{\sum_P x(P)y_{\mathcal{L}_T}(P)}{\operatorname{cross}_T(\alpha)}, \]
 when $ b_{T}(\mathcal{L}_{T}, i_{d})=-1$,
\end{itemize}
where $M$ is a obtained by multiplying together a sequence of matrices of the four forms mentioned forms \eqref{eqn:ccw+1}, \eqref{eqn:ccw-1}, \eqref{eqn:cw+1}, \eqref{eqn:cw-1} and can be interpreted as the matrix $M_d$ for $\overline{\mathcal{G}}_{T, \gamma}$ by comparing with Definition~\ref{definition:tilematrix}. This completes the proof.
\end{proof}

If $\rho$ is a non-standard $M$-path for the one-sided closed curve $\alpha$, then $\rho$ can be deformed into a standard $M$-path $\rho_{\alpha}$ by the local adjustments found in \cite[Lemma 4.8]{MW13} and Lemma \ref{lemma:independenceofpath2}.

\begin{Lemma}\label{lemma:positivecoef}
If we use the standard $M$-path, Definition~{\rm\ref{definition:Mpath2}}, for a one-sided closed curve $\alpha$ with no contractible kinks, then every coefficient of $\chi_{\alpha, \mathcal{L}_T}$ is positive.
\end{Lemma}
\begin{proof}
Given a one-sided closed curve $\alpha$ with no contractible kinks, let $\rho_{\alpha}$ denote the corresponding standard $M$-path. The proof of this lemma comes from the simple observation that all entries of matrices of the forms \eqref{eqn:ccw+1}, \eqref{eqn:ccw-1}, \eqref{eqn:cw+1}, \eqref{eqn:cw-1}, \eqref{eqn:finalstep1} and \eqref{eqn:finalstep2} are positive. As such, $\chi_{\alpha, \mathcal{L}_T}$ is a sum of positive terms and is thus positive.
\end{proof}

An immediate corollary of Theorem~\ref{thm:matrixformula} and Lemma \ref{lemma:positivecoef} is the following.

\begin{Corollary}
The quantity $\chi_{\alpha, \mathcal{L}_T}$ is a Laurent polynomial with all coefficients positive. This verifies and provides another proof of the positivity result by Wilson in {\rm\cite{Wil19}}.
\end{Corollary}

Next, we need to make sure that if we fix a principal lamination $\mathcal{L}_T$, then the isotopic representations of a one-sided curve $\alpha$ share the same Laurent polynomial. That is, if $\alpha$ and~$\alpha'$ are isotopic one-sided closed curves, then $\chi_{\alpha, \mathcal{L}_T}=\chi_{\alpha', \mathcal{L}_T}$. Throughout this segment, we will write \[\chi_{\alpha, \mathcal{L}_T}=\frac{\sum_P x(P)y_{\mathcal{L}_T}(P)}{\operatorname{cross}_T(\alpha)}\qquad \text{ and } \qquad \chi_{\alpha', \mathcal{L}_T}=\frac{\sum_{P} x'(P)y_{\mathcal{L}_T}'(P)}{\operatorname{cross}_T(\alpha')}.\]

\begin{figure}
 \centering
 \includegraphics[scale=.5]{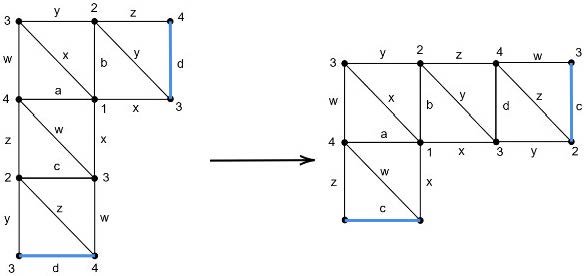}
 \caption{Band graph of $\beta$ obtained by rotating $\alpha$. Despite the band graphs being non-isomorphic, they provide the same expansion due to $\beta$ having a $Z$-intersection with~$z$.}\label{fig:rotationgraph}
\end{figure}

While this statement may seem straightforward, there is some difficulty that occurs due to the fact that, in most situations, the band graphs of $\alpha$ and $\alpha'$ are not isomorphic to one another and hence yield different poset structures. This is illustrated in Figure~\ref{fig:rotationgraph}.

In proving $\chi_{\alpha, \mathcal{L}_T}=\chi_{\alpha', \mathcal{L}_T}$, the first case we will consider is when $\alpha$ and $\alpha'$ are ``reflections'' of one another. Given an arc $\alpha$, its reflection $\alpha'$ crosses the arcs of the triangulation in the same order, but now every $Z$-intersection of $\alpha$ is now an $S$-intersection and every $S$-intersection is now a $Z$-intersection. This has the effect of taking $\alpha$ and rotating it $180^{\circ}$ about the cross cap. For an example of a reflection, compare $\alpha$ in Figure \ref{fig:mobiuspath} with $\alpha'$ in Figure \ref{fig:maxpath}.

\begin{Lemma}\label{lemma:reflections}
 Let $\alpha$ be a one-sided closed curve and $\mathcal{L}_T$ be a principal lamination on the triangulation $T$. Let $\alpha'$ be the reflection of $\alpha$ about the crosscap. Then $\chi_{\alpha, \mathcal{L}_T}=\chi_{\alpha', \mathcal{L}_T}$.
\end{Lemma}
\begin{proof}
 By definition, we have $\operatorname{cross}_T(\alpha)=\operatorname{cross}_T(\alpha')$; furthermore, one can obtain the band graph $\widetilde{\mathcal{G}}_{T, \alpha'}$ of $\alpha'$ by taking $\widetilde{\mathcal{G}}_{T, \alpha}$ and reversing the relative orientation of each tile. This ensures that the set of perfect matchings on both graphs are the same, and $x(P)=x'(P)$ for each perfect matching $P$. Note, the minimal matching on $\widetilde{\mathcal{G}}_{T, \alpha}$ is the maximum matching on $\widetilde{\mathcal{G}}_{T, \alpha'}$ (and vice versa); however, reflecting $\alpha$ across the crosscap has the effect of changing $S$-intersections into $Z$-intersections (and vice versa).

 As a corollary, $y_{\mathcal{L}_T}(P)=y_{\mathcal{L}_T}'(P)$ for this specific matching. The remaining matchings on both band graphs are found by taking $P$, flipping the orientation of certain diagonals and repeating this process until conclusion. We thus have \[\chi_{\alpha, \mathcal{L}_T}=\frac{\sum_P x(P)y_{\mathcal{L}_T}(P)}{\operatorname{cross}_T(\alpha)}=\frac{\sum_{P} x'(P)y_{\mathcal{L}_T}'(P)}{\operatorname{cross}_T(\alpha')}=\chi_{\alpha', \mathcal{L}_T}.\tag*{\qed}\] \renewcommand{\qed}{}
\end{proof}

\begin{Proposition}\label{prop:rotation}
 Let $\alpha$ and $\alpha'$ be isotopic one-sided closed curves. Then $\chi_{\alpha, \mathcal{L}_T}=\chi_{\alpha', \mathcal{L}_T}$.
\end{Proposition}
\begin{proof}
 We begin by setting up some notation. Let $\tau_1, \ldots, \tau_n$ be the arcs of our triangulation $T$ and let $\mathcal{L}_T$ be a principal lamination. As a reminder, each $\tau_i$ lifts to two arcs in the orientable double cover that will have opposite signs with respect to the principal lamination. Let $\overline{\tau_i}$ be the copy of $\tau_i$ that lifts to an $S$-shape intersection and $\widetilde{\tau_i}$ be the copy that lifts to a $Z$-shape intersection.

 To prove this result, we show $\chi$ is invariant under ``rotations'' of the one-sided closed curve around the crosscap, where a ``rotation'' is really just a change of basepoint of the curve $\alpha$. In other words, if $\tau_i$ and $\tau_{i+1}$ are the first two arcs that $\alpha$ intersects and $\tau_{i-1}$ is the final arc~$\alpha$ intersects, then we will assume $\alpha'$ intersects $\tau_{i+1}$ first and $\tau_{i}$ last (see Figure~\ref{fig:rotation}) and show~${\chi_{\alpha, \mathcal{L}_T}=\chi_{\alpha', \mathcal{L}_T}}$. Note that the $\tau_i$'s in $\alpha$ and $\alpha'$ correspond to arcs on the orientable surface with opposite sign with respect to $\mathcal{L}_T$. By Lemma \ref{lemma:reflections}, we may assume without loss of generality, that $\tau_{i+1}$ lies counterclockwise to $\tau_{i}$. This implies that on the other side of the crosscap, $\tau_{i+1}$ lies clockwise to $\tau_{i}$ (otherwise, we can reflect across the crosscap). Further, we will assume the~$\tau_i$ that $\alpha$ crosses lifts to $\overline{\tau_i}$, while the $\tau_i$ that $\alpha'$ crosses lifts to $\widetilde{\tau_i}$.

 By Theorem~\ref{thm:matrixformula} and the fact $\operatorname{tr}(UV)=\operatorname{tr}(VU)$ for two matrices $U$ and $V$, we have
 \begin{align*}
 \chi_{\alpha, \mathcal{L}_T}&=\operatorname{tr}\left(
 \begin{bmatrix}
 0 & x_i\\
 \dfrac{1}{x_i} & 0
 \end{bmatrix}\cdot M\cdot
 \begin{bmatrix}
 1 & 0\\
 \dfrac{x_{a_i}}{x_i x_{i+1}} & 1
 \end{bmatrix}
 \begin{bmatrix}
 1 & 0\\
 0 & y_i
 \end{bmatrix}\right)\\
 &=\operatorname{tr}\left(
 \begin{bmatrix}
 1 & 0\\
 \dfrac{x_{a_i}}{x_i x_{i+1}} & 1
 \end{bmatrix}
 \begin{bmatrix}
 1 & 0\\
 0 & y_i
 \end{bmatrix}
 \begin{bmatrix}
 0 & x_i\\
 \dfrac{1}{x_i} & 0
 \end{bmatrix}\cdot M
 \right)
 = \operatorname{tr}\left(
 \begin{bmatrix}
 0 & x_i\\
 \dfrac{y_i}{x_i} & \dfrac{x_{a_i}}{x_{i+1}}
 \end{bmatrix}\cdot M
 \right)\\
 &=\operatorname{tr}\left(
 \begin{bmatrix}
 0 & x_{i+1}\\
 \dfrac{1}{x_{i+1}} & 1
 \end{bmatrix}
 \begin{bmatrix}
 1 & 0\\
 \dfrac{x_i}{x_{a_i} x_{i+1}} & 1
 \end{bmatrix}
 \begin{bmatrix}
 0 & x_{a_i}\\
 -\dfrac{1}{x_{a_i}} & 0
 \end{bmatrix}
 \begin{bmatrix}
 1 & 0\\
 \dfrac{x_{i+1}}{x_{a_i}x_{i}} & 1
 \end{bmatrix}
 \begin{bmatrix}
 y_i & 0\\
 0 & 1
 \end{bmatrix}\cdot M
 \right)\\
 & = \chi_{\alpha', \mathcal{L}_T},
\end{align*}

\noindent where $M$ is a product of matrices determined by the standard $M$-path of both $\alpha$ and $\alpha'$, i.e., it is a product of matrices \eqref{eqn:ccw+1}, \eqref{eqn:ccw-1}, \eqref{eqn:cw+1}, \eqref{eqn:cw-1}. The above proof holds, even if we were to switch the signs of $\overline{\tau_i}$ and $\widetilde{\tau_i}$ with respect to $\mathcal{L}_T$.
\end{proof}

\begin{Remark}
One could also prove this by making use of the band graphs of these curves. Namely, a rotation of the $\alpha$ amounts to moving the first tile of the graph to the end of the graph and then flipping the sign of the lamination. This method is far more difficult to prove in generality than the matrix method used above.
\end{Remark}

\subsection{M\"obius band example}\label{section:example}

In this subsection, we illustrate the techniques used above
in a few examples involving the M\"obius Band. We will first compute the Laurent expansion
by computing the trace of the $M$-path matrix and then we
will compare it to the Laurent expansion obtained from the
good matching enumerator.

\noindent \textbf{Notation:}
 Throughout this section, if we have an arc $x$ in the triangulation, we let $X$ denote the coefficient $y_x$ corresponding to $x$.

\begin{figure}[!ht]\label{fig:alphafigure}
\centering
\includegraphics[scale=.4]{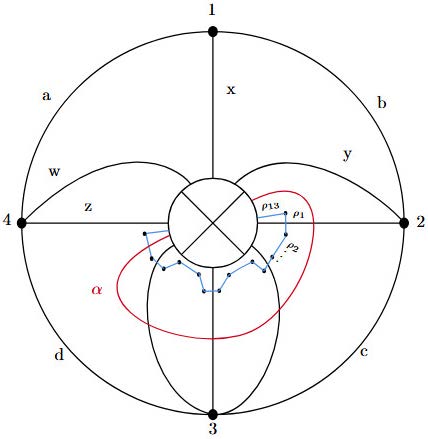}
\caption{The one-sided curve $\alpha$ on $M_4$ along with the standard $M$-path which begins with $\rho_1$ and ends with $\rho_{13}$.}\label{fig:mobiuspath}
\end{figure}

The one-sided curve $\alpha$ on $M_4$ and each step of the
standard $M$-path are shown in Figure~\ref{fig:mobiuspath}. We assume each arc $w$, $x$, $y$, $z$ that $\alpha$ intersects lifts to an $S$-intersection in the orientable double cover and will look at what happens on different laminations later. The matrices associated to each of step of the standard $M$-path are listed here
\begin{alignat*}{4}
& M( \rho _{1}) = \begin{bmatrix}
1 & 0\\
0 & Z
\end{bmatrix} ,\qquad&& M( \rho _{2}) = \begin{bmatrix}
1 & 0\\
\dfrac{c}{wz} & 1
\end{bmatrix} ,\qquad &&M( \rho _{3}) = \begin{bmatrix}
1 & 0\\
0 & W
\end{bmatrix} ,&\\ &M( \rho _{4}) = \begin{bmatrix}
1 & 0\\
\dfrac{x}{aw} & 1
\end{bmatrix} ,
\qquad&&
M( \rho _{5}) = \begin{bmatrix}
0 & a\\
\dfrac{-1}{a} & 0
\end{bmatrix} ,\qquad&& M( \rho _{6}) = \begin{bmatrix}
1 & 0\\
\dfrac{w}{ax} & 1
\end{bmatrix} ,&\\ &M( \rho _{7}) = \begin{bmatrix}
1 & 0\\
0 & X
\end{bmatrix} ,\qquad&& M( \rho _{8}) = \begin{bmatrix}
1 & 0\\
\dfrac{y}{xb} & 1
\end{bmatrix} ,
\qquad&&
M( \rho _{9}) = \begin{bmatrix}
0 & b\\
\dfrac{-1}{b} & 0
\end{bmatrix} ,&\\ &M( \rho _{10}) = \begin{bmatrix}
1 & 0\\
\dfrac{x}{by} & 1
\end{bmatrix} ,\qquad&& M( \rho _{11}) = \begin{bmatrix}
1 & 0\\
0 & Y
\end{bmatrix} ,\qquad&& M( \rho _{12}) = \begin{bmatrix}
1 & 0\\
\dfrac{d}{yz} & 1
\end{bmatrix} ,&
\\
 &M( \rho _{13}) = \begin{bmatrix}
0 & z\\
\dfrac{1}{z} & 0
\end{bmatrix}.&&&&&
\end{alignat*}
Therefore, we have
\begin{align*}
 M( \rho _{\alpha }) ={}&M( \rho _{13}) \cdot [ M( \rho _{12}) \cdot M( \rho _{11})] \cdot [ M( \rho _{10}) \cdot M( \rho _{9}) \cdot M( \rho _{8}) \cdot M( \rho _{7})] \cdot [ M( \rho _{6}) \cdot M( \rho _{5}) \\
&	\cdot	M( \rho _{4}) \cdot M( \rho _{3})] \cdot [ M( \rho _{2}) \cdot M( \rho _{1})]
\\
={}&\begin{bmatrix}
0 & z\\
\dfrac{1}{z} & 0
\end{bmatrix} \cdot \begin{bmatrix}
1 & 0\\
\dfrac{d}{yz} & Y
\end{bmatrix} \cdot \begin{bmatrix}
\dfrac{y}{x} & bX\\
0 & \dfrac{xX}{y}
\end{bmatrix} \cdot \begin{bmatrix}
\dfrac{x}{w} & aW\\
0 & \dfrac{wW}{x}
\end{bmatrix} \cdot \begin{bmatrix}
1 & 0\\
\dfrac{c}{wz} & Z
\end{bmatrix}
\\
={}&\begin{bmatrix}
\dfrac{\begin{array}{@{}c@{}} cwxzWXY+bcdwWX+bcdwWX\\ +acdyW+dxyz\end{array}}{wxyz} & \dfrac{wxzWXYZ+bdwWXZ+adyWZ}{xy}\vspace{2mm}\\
\dfrac{bcwWX+acyW+xyz}{wxz^{2}} & \dfrac{bwWXZ+ayWZ}{xz}
\end{bmatrix}
\\
\Rightarrow{}& \mathbf{\text{Trace}}( M( \rho _{\alpha }))=\chi_{\alpha, \mathcal{L}_T} \\
={}&\frac{cwxzWXY+bw^{2} yWXZ+bcdwWX+awy^{2} WZ+acdyW+dxyz}{wxyz},
\end{align*}
where the closed brackets, in the $M(\rho_{\alpha})$ product, group the matrices by
triangle and the final matrix, $M(\rho_{13})$, corresponds to going through the crosscap.

\begin{figure}[!ht]
\centering
\includegraphics[scale=.6]{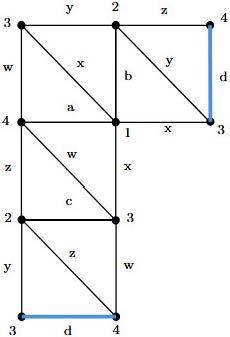}
\vspace{-3mm}

\caption{Band graph $\widetilde{\mathcal{G}}_{\alpha, T}$ for $M_4$.}\label{fig:snakegraph}
\end{figure}

\begin{figure}[!ht]
\centering
\includegraphics[scale=.30]{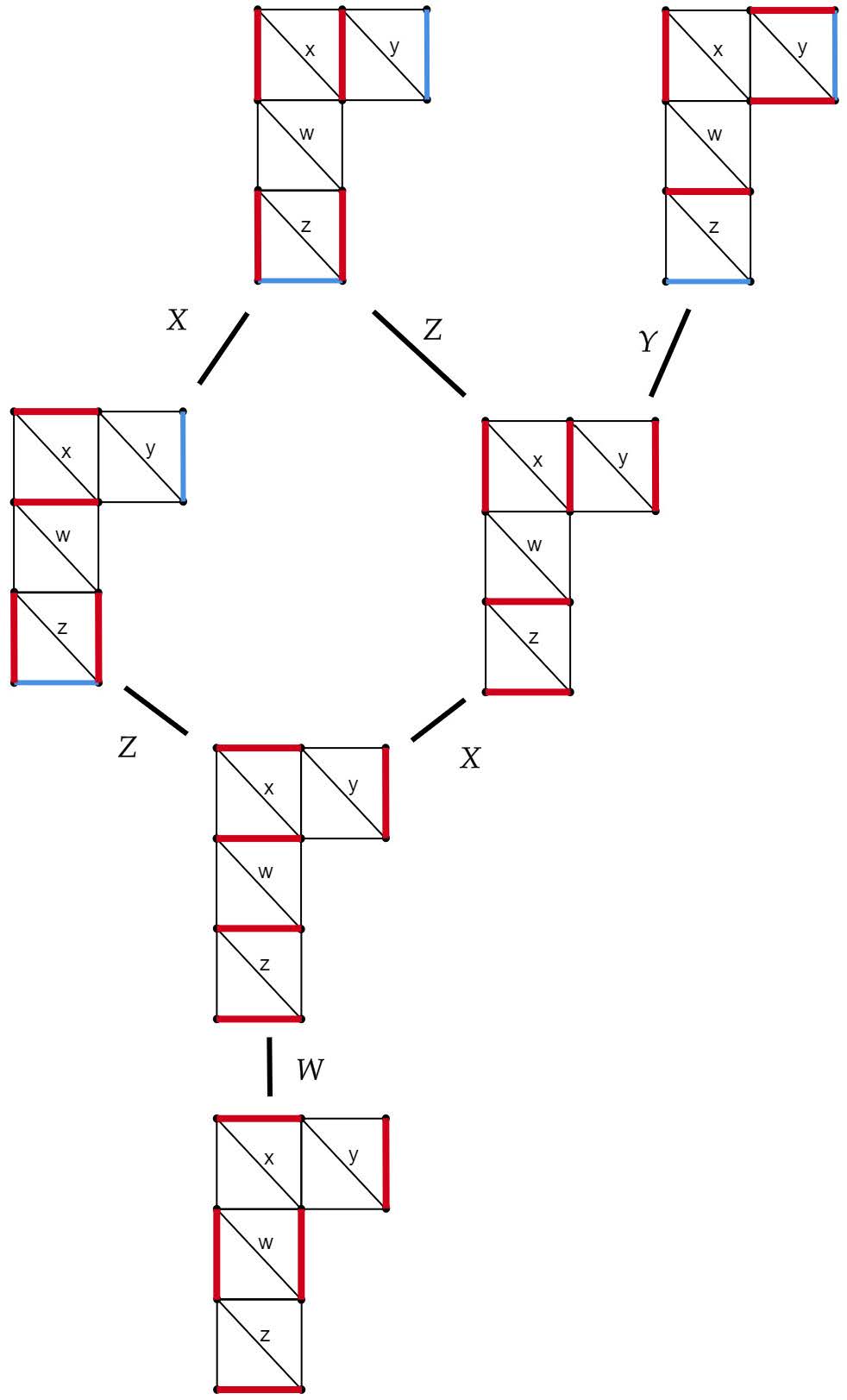}

\caption{Poset structure associated with $\overline{\mathcal{G}}_{\alpha, T}$.}\label{fig:posetgraph}
\end{figure}

Meanwhile, we have the band graph $\overline{\mathcal{G}}_{\alpha, T}$ shown in Figure~\ref{fig:snakegraph}. With this band graph, the good matchings are demonstrated in Figure~\ref{fig:posetgraph}.

 Computing the good matching enumerator, we find
\[\frac{\sum_{P} x(P)y_{\mathcal{L}_T}(P)}{\operatorname{cross}_T(\alpha)}
=\frac{cwxzWXY+bw^{2} yWXZ+bcdwWX+awy^{2} WZ+acdyW+dxyz}{wxyz},\]
and we see the two methods of computing the expansion agree with one another.

\subsection{Reflection and rotation example}

\begin{figure}[!ht]\label{fig:alpha'fig}
 \centering
 \includegraphics[scale=.62]{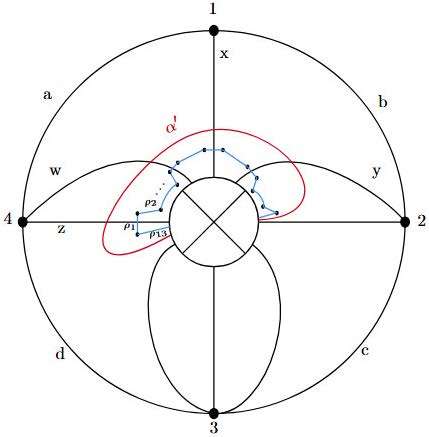}
 \caption{$M$-path after reflecting $\alpha$ from the previous example across the crosscap}\label{fig:maxpath}
\end{figure}

If we instead considered the curve given by reflecting $\alpha$ across the crosscap, which we will denote by $\alpha'$, then we get the $M$-path found in Figure~\ref{fig:maxpath}. The expansion is given by
\begin{align*}
\chi_{\alpha', \mathcal{L}_T}&=\operatorname{tr}\left(
\begin{bmatrix}
 0 && z\\
 \dfrac{1}{z} && 0
\end{bmatrix}
\begin{bmatrix}
 \dfrac{zY}{y} && d\\
 0 && \dfrac{y}{z}
\end{bmatrix}
\begin{bmatrix}
 X && 0\\
 \dfrac{bX}{xy} && 1
\end{bmatrix}
\begin{bmatrix}
 W && 0\\
 \dfrac{aW}{wx} && 1
\end{bmatrix}
\begin{bmatrix}
 \dfrac{wZ}{z} && c\\
 0 && \dfrac{z}{w}
\end{bmatrix}\right)\\
&=\frac{cwxzWXY+bw^{2} yWXZ+bcdwWX+awy^{2} WZ+acdyW+dxyz}{wxyz},
\end{align*}
which agrees with out previous computation. Observe that when flipping across the crosscap, every counterclockwise sequence becomes clockwise and vice-versa. Additionally, the sign of the principal lamination changes as we are now intersecting arcs that lift to $Z$-intersections, so we must change each matrix of type two in the standard $M$-path.

\begin{figure}[!ht]
 \centering
 \includegraphics[scale=.7]{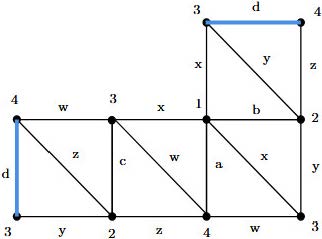}
 \caption{Band graph $\widetilde{\mathcal{G}}_{\alpha', T}$.}
 \label{fig:maximalband}
\end{figure}

The corresponding band graph is shown in Figure~\ref{fig:maximalband}. Observe that it is the same as the band graph shown in Figure~\ref{fig:snakegraph} with the exception that we have used the opposite relative orientation.

\begin{figure}[!ht]
 \centering
 \includegraphics[scale=.55]{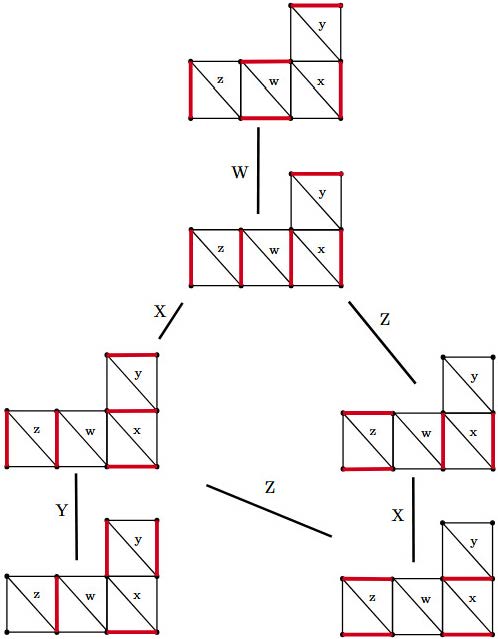}
 \caption{Perfect matchings on $\widetilde{\mathcal{G}}_{\alpha', T}$ which creates a different poset structure than $\widetilde{\mathcal{G}}_{\alpha, T}$ using the flip definition from \cite{MSW11}.} \label{fig:posetforalpha'}
\end{figure}

\begin{figure}[!ht]
 \centering
 \includegraphics[scale=.6]{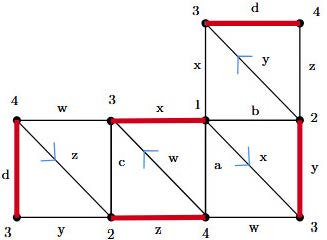}
 \caption{Maximal Matching on $\widetilde{\mathcal{G}}_{\alpha', T}$.}
 \label{fig:maximalmatch}
\end{figure}

Observe that the minimal matching shown in Figure~\ref{fig:posetgraph} (the bottommost graph of the figure), becomes the maximal matching for $\widetilde{\mathcal{G}}_{\alpha', T}$ (see Figure~\ref{fig:maximalmatch}); however, as we are now lifting to $Z$-intersections, the coefficient monomial remains the exact same as the previous example. A~similar logic applies to the other good matchings on $\widetilde{\mathcal{G}}_{\alpha', T}$ and the set of good matchings is shown to give comparison in Figure~\ref{fig:posetforalpha'}.

\begin{figure}[!ht]
 \centering
 \includegraphics[scale=.57]{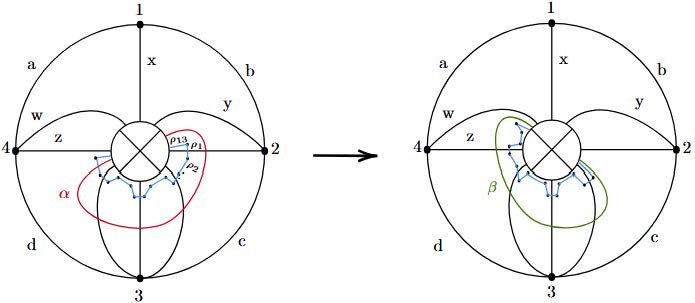}
 \caption{Rotation of $\alpha$ to an isotopic one-sided closed curve $\beta$.}\label{fig:rotation}
\end{figure}

Next, we will look at a ``rotation'' of $\alpha$. Let $\alpha$ be the arc depicted in Figure~\ref{fig:mobiuspath} and let $\beta$ be the one-sided closed curve shown in Figure~\ref{fig:rotation}. The matrices corresponding to the individual steps of the standard $M$-path for $\beta$ are listed below. The first ten matrices appear in the matrix product for $\alpha$ as well,
\begin{alignat*}{4}
&M( \rho_{1}) = \begin{bmatrix}
1 & 0\\
0 & W
\end{bmatrix} ,\qquad&&
 M( \rho _{2}) = \begin{bmatrix}
1 & 0\\
\dfrac{x}{aw} & 1
\end{bmatrix} ,\qquad&& M( \rho _{3}) = \begin{bmatrix}
0 & a\\
\dfrac{-1}{a} & 0
\end{bmatrix} ,&\\ &M( \rho _{4}) = \begin{bmatrix}
1 & 0\\
\dfrac{w}{ax} & 1
\end{bmatrix} ,\qquad&&
M( \rho _{5}) = \begin{bmatrix}
1 & 0\\
0 & X
\end{bmatrix} ,\qquad&& M( \rho _{6}) = \begin{bmatrix}
1 & 0\\
\dfrac{y}{bx} & 1
\end{bmatrix} ,&\\ &M( \rho _{7}) = \begin{bmatrix}
0 & b\\
\dfrac{-1}{b} & 0
\end{bmatrix} ,\qquad&& M( \rho _{8}) = \begin{bmatrix}
1 & 0\\
\dfrac{x}{by} & 1
\end{bmatrix} ,
\qquad&&
M( \rho _{9}) = \begin{bmatrix}
1 & 0\\
0 & Y
\end{bmatrix} ,&\\ &M( \rho _{10}) = \begin{bmatrix}
1 & 0\\
\dfrac{d}{yz} & 1
\end{bmatrix} ,\qquad&& M( \rho _{11}) = \begin{bmatrix}
Z & 0\\
0 & 1
\end{bmatrix} ,\qquad&& M( \rho _{12}) = \begin{bmatrix}
1 & 0\\
\dfrac{w}{cz} & 1
\end{bmatrix},&
\\
&M( \rho _{13}) = \begin{bmatrix}
0 & c\\
\dfrac{-1}{c} & 0
\end{bmatrix} ,\qquad&& M( \rho _{14}) = \begin{bmatrix}
1 & 0\\
\dfrac{z}{cw} & 1
\end{bmatrix} ,\qquad&& M( \rho _{15}) = \begin{bmatrix}
0 & w\\
\dfrac{1}{w} & 0
\end{bmatrix}.&
\end{alignat*}

 The proof of Proposition~\ref{prop:rotation} shows that the trace of the matrix product associated to $\rho_{\alpha}$ agrees with the trace of the matrix product associated to $\rho_{\beta}$. This can be seen explicitly by computing the matrix product via \texttt{Macaulay2}. One finds
\begin{align*}
 \text{Trace}(M(\rho_{\beta}))&=\text{Trace}(M(\rho_{15})\cdots M(\rho_1))\\
 &=\begin{bmatrix}
 \dfrac{d}{w} && \dfrac{wxzWXY+bdwWX+adyW}{xy}\vspace{2mm}\\
 \dfrac{wyZ+cd}{w^2z} && \dfrac{\begin{array}{@{}c@{}}cwxzWXY+bw^2yWXZ+bcdwWX\\ +awy^2WZ+acdyW\end{array}}{wxyz}
 \end{bmatrix}\\
 & = \frac{cwxzWXY+bw^{2} yWXZ+bcdwWX+awy^{2} WZ+acdyW+dxyz}{wxyz}.
\end{align*}
A similar computation can be done with the other rotations of $\alpha$.

\section{Skein relations for one-sided closed curves}\label{section:skeinrelations}

To conclude this paper, we briefly show how the theory we have developed throughout this paper can be used to prove skein relations on non-orientable surfaces. Skein relations on non-orientable surfaces have already been proven for $x$-coordinates in \cite[Section 4.3]{DP15}; however, with our matrix formulae, we can show these equations hold true when considering coefficients that come from principal laminations.

There are only a few cases that we need to consider, and most follow from basic adjustments to the arguments found in \cite{MW13} using Sections \ref{section:matrixorientable} and \ref{section:expansionformulanonorientable} of this paper. Throughout this section, we will let $e(\gamma, \tau)$ denote the number of crossings between the generalized arc/loop/one-sided closed curve $\gamma$ and the arc $\tau$.

\subsection*{Case one: Intersection of generalized arcs} In this case, we do not consider one-sided closed curves, so we can simply lift to the orientable double cover and apply \cite[Proposition 6.4]{MW13}. Considering general principal laminations rather than principal coefficients, the statement becomes the following.

\begin{Proposition}
 Let $\gamma_1$ and $\gamma_2$ be two generalized arcs which intersect each other at least once; let $x$ be a point of intersection; and let $\alpha_1$, $\alpha_2$ and $\beta_1$, $\beta_2$ be the two pairs of arcs obtained by resolving the intersection of $\gamma_1$ and $\gamma_2$ at $x$. Then
 \[\overline{\chi}_{\gamma_1, \mathcal{L}_T}\overline{\chi}_{\gamma_2, \mathcal{L}_T}=\pm \overline{\chi}_{\alpha_1, \mathcal{L}_T}\overline{\chi}_{\alpha_2, \mathcal{L}_T}\prod_{i=1}^n y_i^{b_T(\mathcal{L}_T, \tau_i)(c_i-a_i)/2}\pm \overline{\chi}_{\beta_1, \mathcal{L}_T}\overline{\chi}_{\beta_2, \mathcal{L}_T}\prod_{i=1}^n y_i^{b_T(\mathcal{L}_T, \tau_i)(c_i-b_i)/2},\]
 where
 $c_i=e(\gamma_1, L_i)+e(\gamma_2, L_i), a_i=e(\alpha_1, L_i)+e(\alpha_2, L_i)$ and $b_i=e(\beta_1, L_i)+e(\beta_2, L_i)$.
\end{Proposition}

\subsection*{Case two: Intersection of generalized arc/loop with generalized \\ loop/one-sided closed curve}

When $\gamma_1$ is a generalized arc or loop and $\gamma_2$ is a generalized loop, then \cite[Proposition 6.5]{MW13} applies after lifting to the orientable double cover. When $\gamma_2$ is a one-sided closed curve, we note that the only new type of step in the standard $M$-path is the elementary step of type 3$'$, from Definition~\ref{definition:newstep}, when going through the crosscap at the end, see Definition~\ref{definition:Mpath2}. Specifically, the steps of type~1 and~2 remain the same as in the orientable case, meaning we can directly apply the arguments from \cite[Lemma~6.10 and Proposition~6.5]{MW13} to prove Proposition~\ref{prop:skeincurvewitharc}.

\begin{figure}[!ht]
 \centering
 \includegraphics[scale=.55]{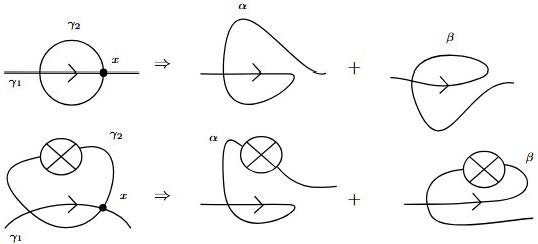}
 \caption{Examples of intersections with curves and their resolutions.}\label{fig:Skein Relation: Intersection of Curve with Arc}
\end{figure}

\begin{Proposition}\label{prop:skeincurvewitharc}
 Let $\gamma_1$ be a generalized arc or loop and let $\gamma_2$ be a generalized loop or one-sided closed curve, such that $\gamma_1$ and $\gamma_2$ intersect each other at least once; let $x$ be a point of intersection; and let $\alpha$ and $\beta$ be the two arcs/curves obtained by resolving the intersection of $\gamma_1$ and $\gamma_2$ at $x$. Then \[\overline{\chi}_{\gamma_1, \mathcal{L}_T}\overline{\chi}_{\gamma_2, \mathcal{L}_T}=\pm \overline{\chi}_{\alpha, \mathcal{L}_T}\prod_{i=1}^n y_i^{b_T(\mathcal{L}_T, \tau_i)(c_i-a_i)/2}\pm \overline{\chi}_{\beta, \mathcal{L}_T}\prod_{i=1}^n y_i^{b_T(\mathcal{L}_T, \tau_i)(c_i-b_i)/2},\]
 where $c_i=e(\gamma_1, L_i)+e(\gamma_2, L_i), a_i=e(\alpha, L_i)$ and $b_i=e(\beta, L_i)$.
\end{Proposition}

\subsection*{Case three: Non-simple generalized arc/loop/one-sided closed curve}

Just like the previous case, when $\gamma$ is a generalized arc or closed curve with a self-intersection at $x$, me may lift to the orientable surface and apply \cite[Proposition~6.6]{MW13}. When $\gamma$ is a one-sided closed curve, using the standard $M$-path for one-sided closed curves, we may directly apply the arguments from \cite[Lemma 6.10 and Proposition 6.6]{MW13} to conclude Proposition~\ref{prop:selfintersection}.

\begin{figure}[!ht]
 \centering
 \includegraphics[scale=.5]{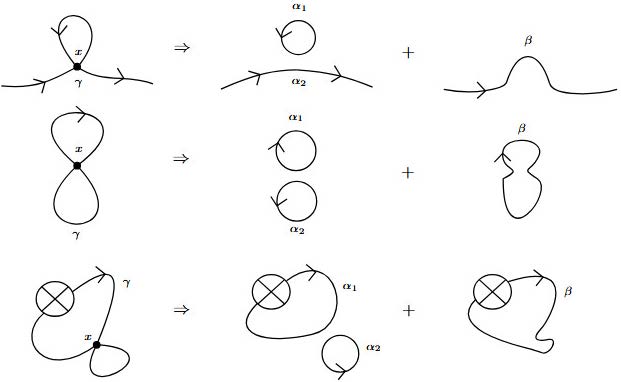}
 \caption{Examples of self-intersection resolution.}\label{fig:Skein Relation: Self-intersection}
\end{figure}

\begin{Proposition}\label{prop:selfintersection}
 Let $\gamma$ be a generalized arc, closed curve or one-sided closed curve with a~self-intersection at $x$. Let $\alpha_1$, $\alpha_2$ and $\beta$ be the generalized arcs/loops/one-sided curves obtained by resolving the intersection at $x$. Then \[\overline{\chi}_{\gamma, \mathcal{L}_T}=\pm \overline{\chi}_{\alpha_1, \mathcal{L}_T}\overline{\chi}_{\alpha_2, \mathcal{L}_T}\prod_{i=1}^n y_i^{b_T(\mathcal{L}_T, \tau_i)(c_i-a_i)/2}\pm \overline{\chi}_{\beta, \mathcal{L}_T}\prod_{i=1}^n y_i^{b_T(\mathcal{L}_T, \tau_i)(c_i-b_i)/2},
 \] where
 $c_i=e(\gamma, L_i)$, $a_i=e(\alpha_1, L_i)+e(\alpha_2, L_i)$ and $b_i=e(\beta, L_i)$.
\end{Proposition}

\subsection*{Case four: Intersection of homotopic one-sided closed curves}

The case that requires the most care occurs when we have two homotopic one-sided closed curves that intersect one another. It's worth noting that any two curves homotopic to a one-sided closed curve will have at least one intersection point, meanwhile, in the orientable case, two homotopic two-sided curves can always be adjusted so that they are disjoint.

For the next proposition, we follow the notation from \cite[Proposition~4.6]{DP15}. Namely, if $\alpha$ is a one-sided closed curve, then $\alpha^2$ will denote the one-sided closed curve of multiplicity two, see Figure~\ref{fig:squareskein} below. In terms of the orientable double cover, one can interpret $\alpha^2$ as the concatenation of the two lifts of $\alpha$ on the orientable surface, as such, $\alpha^2$ is a two-sided closed curve enclosing the crosscap. The following proposition generalizes \cite[Proposition~4.6]{DP15}.

\begin{figure}[!ht]
 \centering
 \includegraphics[scale=.5]{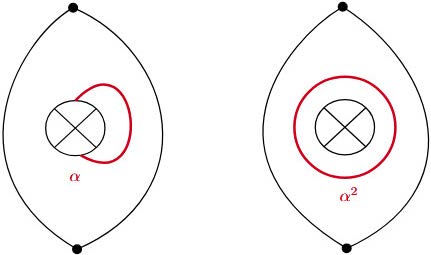}
 \caption{Depiction of $\alpha^2$ on the M\"obius band with two marked points.}
 \label{fig:squareskein}
\end{figure}

\begin{Proposition}
 Let $\alpha$ be a one-sided closed curve, and let $\alpha^2$ be the two sided-closed curve enclosing the crosscap obtained from resolving an intersection of $\alpha$ with another one-sided closed curve homotopic to $\alpha$, then
$(\overline{\chi}_{\alpha, _T})^2=\overline{\chi}_{\alpha^2, T}-2$.
\end{Proposition}
\begin{proof}
Let $\alpha'$ be the reflection of $\alpha$ about the crosscap, as in Lemma \ref{lemma:reflections}. We can write \[\overline{\chi}_{\alpha, T}=\operatorname{tr}\left(
\begin{bmatrix}
 0 & x_{i_1}\\
 \dfrac{1}{x_{i_1}} & 0
\end{bmatrix}
\cdot \overline{M}
\right),\] where $\overline{M}$ is the reduced standard $M$-path of $\alpha$ with the final step excluded (so we do not go through the crosscap). Similarly, we can write \[\overline{\chi}_{\alpha', T}=\operatorname{tr}\left(
\begin{bmatrix}
 0 & x_{i_1}\\
 \dfrac{1}{x_{i_1}} & 0
\end{bmatrix}
\cdot \overline{M'}
\right).\]

With this notation, we have $\overline{\chi}_{\alpha^2, T}=\operatorname{tr}\bigl(\overline{M'}\cdot \overline{M}\bigr)$. With all of this information in mind, from earlier results from this paper, as well as other basic considerations, we have
\begin{itemize}\itemsep=0pt
 \item $\overline{\chi}_{\alpha, T}=\overline{\chi}_{\alpha', T}$,
 \item $\det\bigl(\overline{M}\bigr)=\det\bigl(\overline{M'}\bigr)=1$,
 \item
 \[
 \det\left(\begin{bmatrix}
 0 & x_{i_1}\\
 \dfrac{1}{x_{i_1}} & 0
\end{bmatrix}\right)=-1\qquad \text{and}\qquad \begin{bmatrix}
 0 & x_{i_1}\\
 \dfrac{1}{x_{i_1}} & 0
\end{bmatrix}^2=I_2\]
\item \[
\begin{bmatrix}
 0 & x_{i_1}\\
 \dfrac{1}{x_{i_1}} & 0
\end{bmatrix}\cdot \overline{M'}\cdot
\begin{bmatrix}
 0 & x_{i_1}\\
 \dfrac{1}{x_{i_1}} & 0
\end{bmatrix}=\overline{M}.\]
\end{itemize}
 Lastly, when $A$, $B$ are $2\times 2$ matrices with $|{\det(A)}|=|{\det(B)}|=1$ and $\operatorname{tr}(A), \operatorname{tr}(B)>0$, we have $\operatorname{tr}(A)\operatorname{tr}(B)=\operatorname{tr}(AB)+\det(B)\operatorname{tr}\bigl(AB^{-1}\bigr)$.

 Putting everything together, we have
 \begin{align*}
 (\overline{\chi}_{\alpha, T})^2 &=\overline{\chi}_{\alpha, T}\cdot \overline{\chi}_{\alpha, T}=\overline{\chi}_{\alpha', T}\cdot \overline{\chi}_{\alpha, T}\\
={}&\operatorname{tr}\left(
\begin{bmatrix}
 0 & x_{i_1}\\
 \dfrac{1}{x_{i_1}} & 0
\end{bmatrix}
\cdot \overline{M'}\right)\cdot \operatorname{tr}\left(
\begin{bmatrix}
 0 & x_{i_1}\\
 \dfrac{1}{x_{i_1}} & 0
\end{bmatrix}
\cdot \overline{M}\right)\\
={}&\operatorname{tr}\left(
\begin{bmatrix}
 0 & x_{i_1}\\
 \dfrac{1}{x_{i_1}} & 0
\end{bmatrix}
\cdot \overline{M'}\right)\cdot \operatorname{tr}\left(\overline{M}\cdot
\begin{bmatrix}
 0 & x_{i_1}\\
 \dfrac{1}{x_{i_1}} & 0
\end{bmatrix}\right)\\
={}&\operatorname{tr}\left(
\begin{bmatrix}
 0 & x_{i_1}\\
 \dfrac{1}{x_{i_1}} & 0
\end{bmatrix}
\cdot \overline{M'}\cdot\overline{M}\cdot
\begin{bmatrix}
 0 & x_{i_1}\\
 \dfrac{1}{x_{i_1}} & 0
\end{bmatrix}\right)\\&-\operatorname{tr}\left(
\begin{bmatrix}
 0 & x_{i_1}\\
 \dfrac{1}{x_{i_1}} & 0
\end{bmatrix}
\cdot \overline{M'}\cdot \left(\overline{M}\cdot \begin{bmatrix}
 0 & x_{i_1}\\
 \dfrac{1}{x_{i_1}} & 0
\end{bmatrix}
\right)^{-1}\right)\\
={}&\operatorname{tr}\left(\begin{bmatrix}
 0 & x_{i_1}\\
 \dfrac{1}{x_{i_1}} & 0
\end{bmatrix}\cdot
\begin{bmatrix}
 0 & x_{i_1}\\
 \dfrac{1}{x_{i_1}} & 0
\end{bmatrix}
\cdot \overline{M'}\cdot\overline{M}\right)\\&-\operatorname{tr}\left(
\begin{bmatrix}
 0 & x_{i_1}\\
 \dfrac{1}{x_{i_1}} & 0
\end{bmatrix}
\cdot \overline{M'}\cdot \begin{bmatrix}
 0 & x_{i_1}\\
 \dfrac{1}{x_{i_1}} & 0
\end{bmatrix}\cdot \overline{M}
^{-1}\right)\\
={}&\operatorname{tr}(\overline{M'}\cdot\overline{M})-\operatorname{tr}\bigl(\overline{M}\cdot\overline{M}^{-1}\bigr)
 =\overline{\chi}_{\alpha^2, T}-\operatorname{tr}(I_2)
 =\overline{\chi}_{\alpha^2, T}-2.\tag*{\qed}
 \end{align*}\renewcommand{\qed}{}
\end{proof}

\appendix

\section{Appendix}
In this appendix, we explicitly demonstrate the skein relation manipulations referred to in Section~\ref{section:quasiclusteralgebras}. We begin by demonstrating the missing step in the computation for the quasi-mutation depicted in the body of the paper. Namely, we show on the double cover why we can push a~loop through the crosscap.

\begin{figure}[!ht]
 \centering
 \includegraphics[width=.8\textwidth]{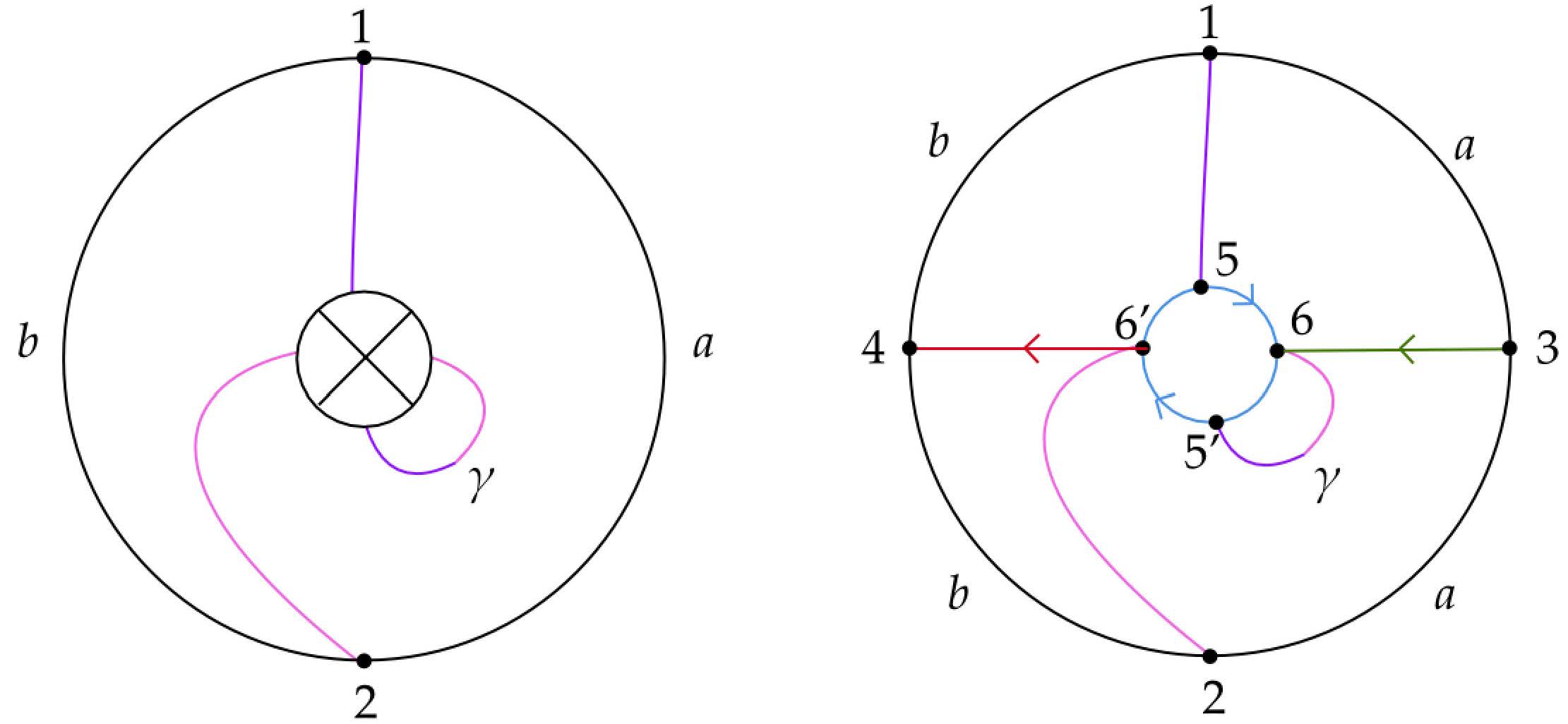}

 \vspace{2mm}

 \centering
 \includegraphics[width=0.9\textwidth]{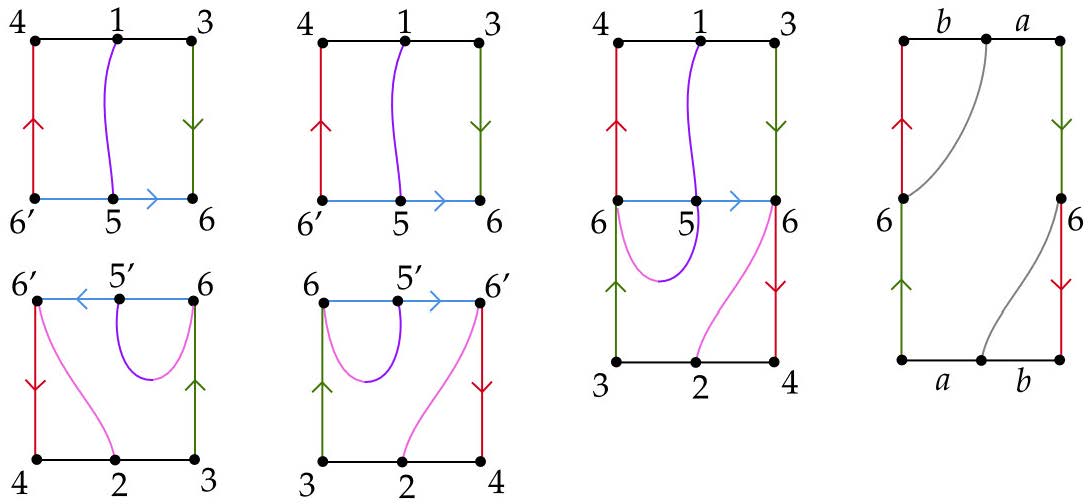}
\end{figure}

Included below is the skein relation derivation of the fourth quasi-mutation relation given in Definition~\ref{definition:quasimutation}. We restate the mutation relation for reference.

Given $t\in T$, the \textit{quasi-mutation} of $\Sigma$ in the direction $T$ as the pair $\mu_t(T,\mathbf{x})=(T',\mathbf{x'})$ where~${T'=\mu_t(T)=T\setminus \{t\}\sqcup\{t'\}}$ and $\mathbf{x'}=\{x_v \mid v\in T'\}$ such that $x_{t'}$ is defined as follows:
If $t$ is an arc separating a triangle with sides $(a,b,t)$ and an annuli with boundary $t$ and one-sided curve~$d$, then the relation is $x_tx_{t'}=(x_a+x_b)^2+x_d^2x_ax_b$.
 \begin{center}
 \includegraphics[width=.5\textwidth]{Figures/skein-reln-4.jpg}
 \end{center}

\begin{figure}[!ht]
 \centering
 \includegraphics[width=.77\textwidth]{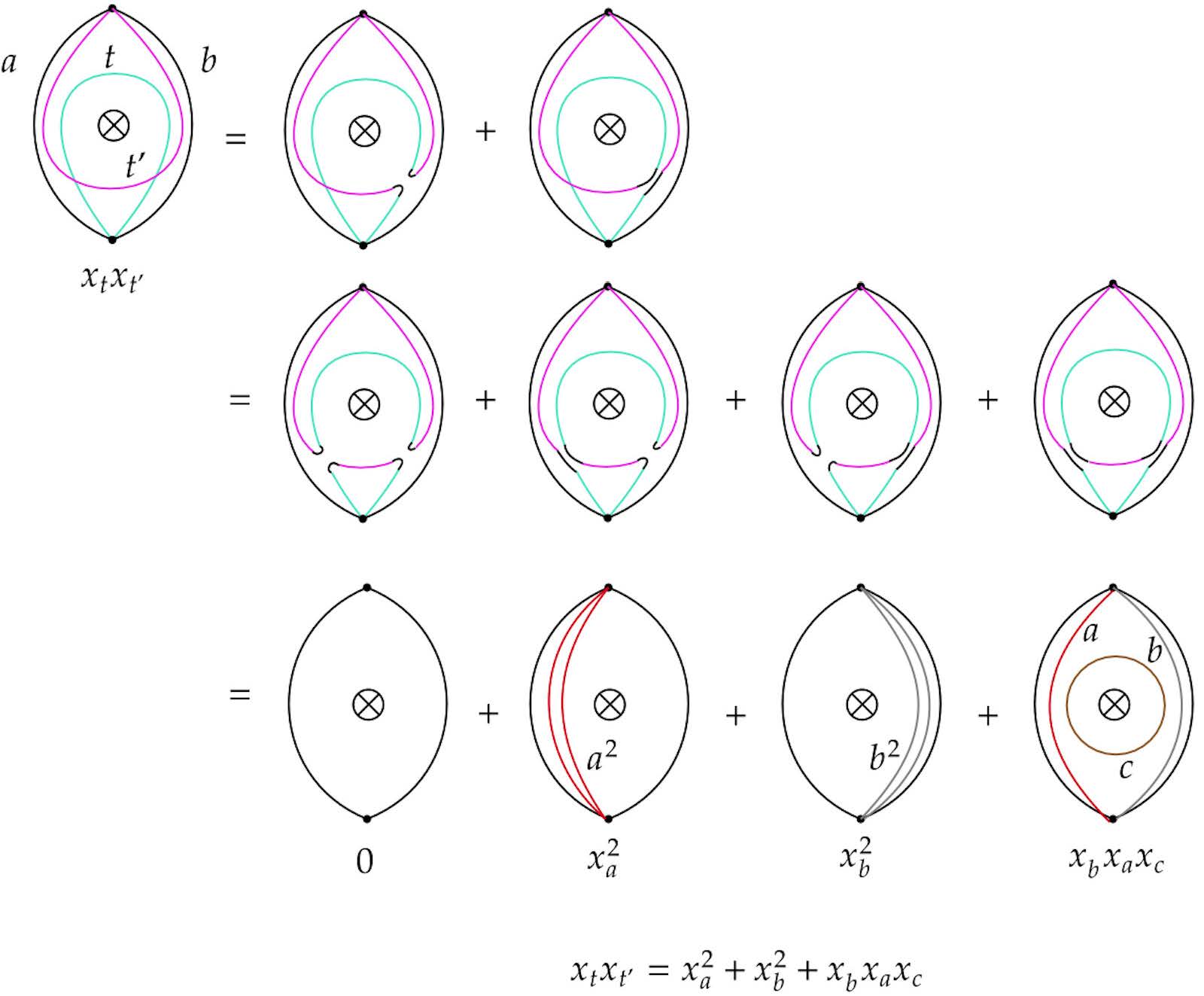}
\end{figure}

The above computation shows that this crossing resolution yields the relation $x_tx_t'=x_a^2+x_b^2+x_bx_ax_c$ where $c$ is a two sided-closed curve enclosing the crosscap. It remains to show that~${x_c=2+x_d^2}$. Let $e$ be the arc passing through the crosscap. We now show the resolution of the crossing of the arc $e$ and the two-sided closed curve $c$ which encloses the crosscap, which yields the relation $x_cx_e=2x_e+x_ax_d+x_bx_d$.

\begin{figure}[!ht]
 \centering
 \includegraphics[width=.8\textwidth]{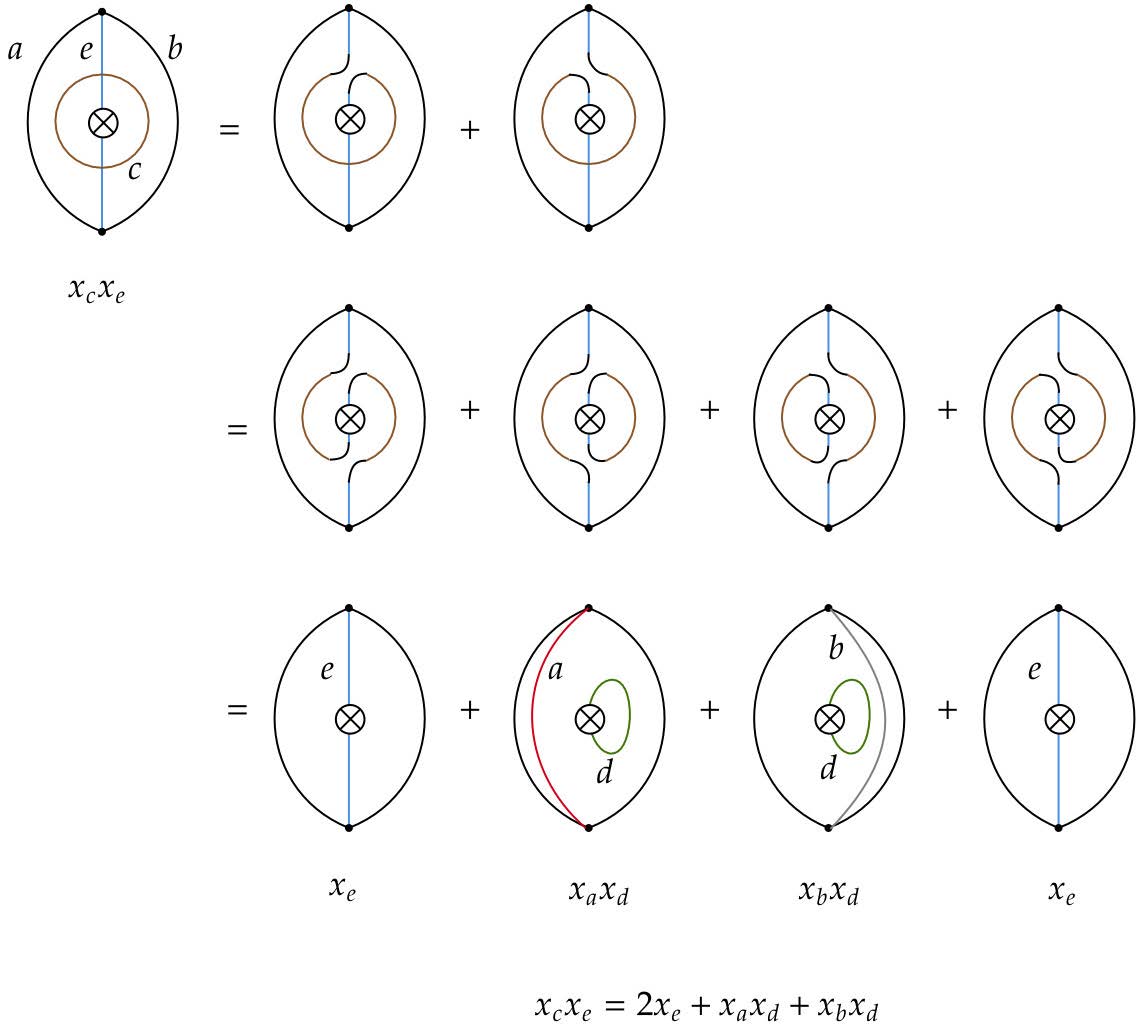}
\end{figure}

Finally, we show that $x_ex_d^2=x_ax_d+x_bx_d$ by resolving the remaining crossings. With this computation, we have $x_cx_e=2x_e+x_ex_d^2$ implying that $x_c=2+x_d^2$. This completes the derivation of the fourth mutation relation from Definition~\ref{definition:quasimutation}.

\begin{figure}[!ht]
 \centering
 \includegraphics[width=.865\textwidth]{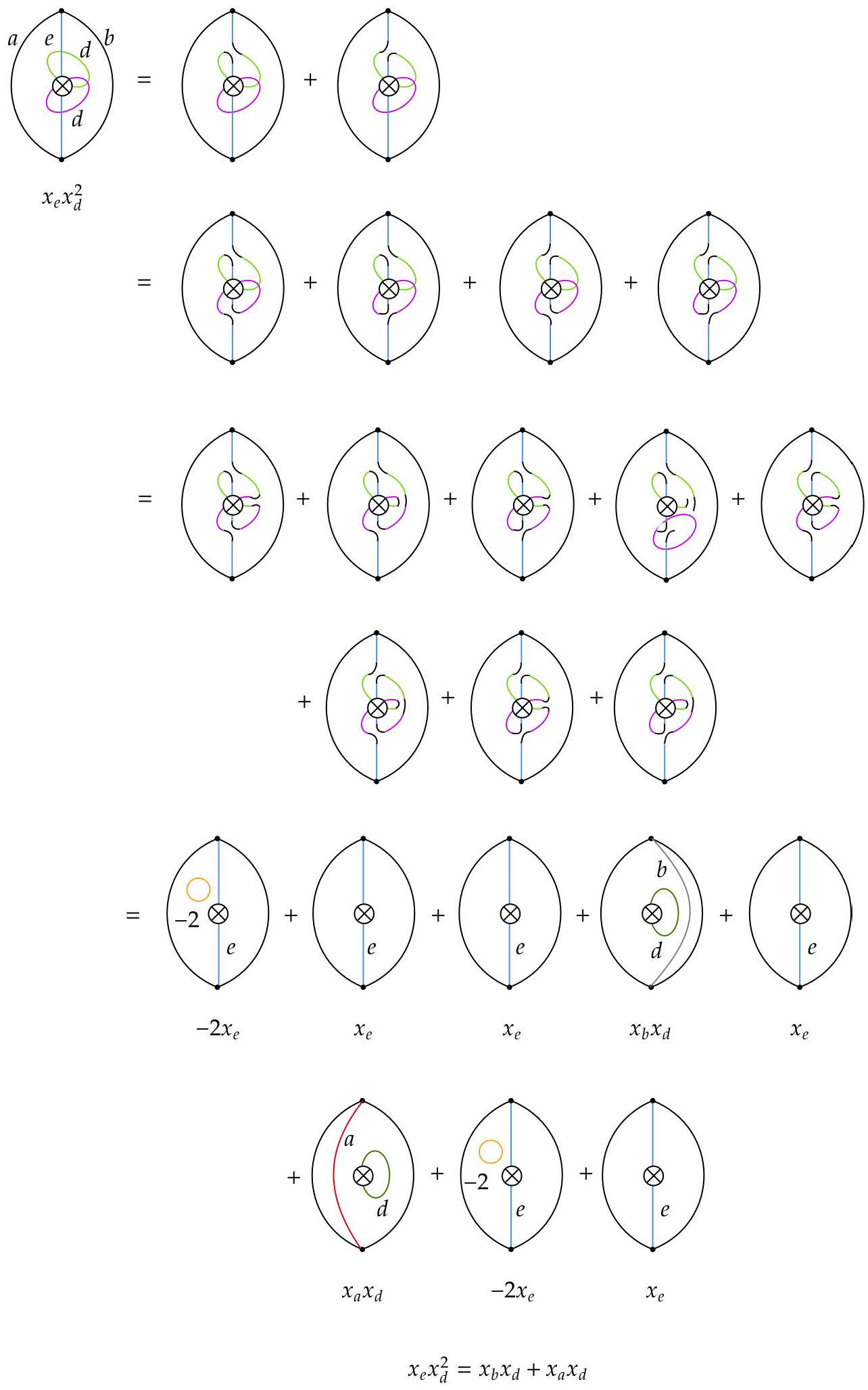}
 \label{fig:quasi-mutation-4c}
\end{figure}

\subsection*{Acknowledgements}

The authors would like to thank Chris Fraser for suggesting this project and helping us with the beginning stages of learning the relevant background. We would also like to thank Gregg Musiker for helpful comments and discussions. Lastly, we would like to thank the anonymous referees during the review process for their useful insights, which led to a clearer, more coherent work.

\pdfbookmark[1]{References}{ref}
\LastPageEnding


\begin{thebibliography}{99}
\footnotesize\itemsep=0pt

\bibitem{DP15}
Dupont G., Palesi F., Quasi-cluster algebras from non-orientable surfaces,
 \href{https://doi.org/10.1007/s10801-015-0586-1}{\textit{J.~Algebraic
 Combin.}} \textbf{42} (2015), 429--472,
 \href{http://arxiv.org/abs/1105.1560}{arXiv:1105.1560}.

\bibitem{FG06}
Fock V., Goncharov A., Moduli spaces of local systems and higher
 {T}eichm\"uller theory,
 \href{https://doi.org/10.1007/s10240-006-0039-4}{\textit{Publ. Math. Inst.
 Hautes \'{E}tudes Sci.}} \textbf{103} (2006), 1--211,
 \href{http://arxiv.org/abs/math.AG/0311149}{arXiv:math.AG/0311149}.

\bibitem{FG07}
Fock V., Goncharov A., Dual {T}eichm\"uller and lamination spaces, in Handbook
 of {T}eichm\"uller {T}heory. {V}ol.~{I}, \textit{IRMA Lect. Math. Theor.
 Phys.}, Vol.~11, \href{https://doi.org/10.4171/029-1/16}{European
 Mathematical Society}, Z\"urich, 2007, 647--684,
 \href{http://arxiv.org/abs/math.DG/0510312}{arXiv:math.DG/0510312}.

\bibitem{FST08}
Fomin S., Shapiro M., Thurston D., Cluster algebras and triangulated
 surfaces.~{I}. {C}luster complexes,
 \href{https://doi.org/10.1007/s11511-008-0030-7}{\textit{Acta Math.}}
 \textbf{201} (2008), 83--146,
 \href{http://arxiv.org/abs/math.RA/0608367}{arXiv:math.RA/0608367}.

\bibitem{FZ02}
Fomin S., Zelevinsky A., Cluster algebras.~{I}. {F}oundations,
 \href{https://doi.org/10.1090/S0894-0347-01-00385-X}{\textit{J.~Amer. Math.
 Soc.}} \textbf{15} (2002), 497--529,
 \href{http://arxiv.org/abs/math.RT/0104151}{arXiv:math.RT/0104151}.

\bibitem{GSV05}
Gekhtman M., Shapiro M., Vainshtein A., Cluster algebras and
 {W}eil--{P}etersson forms,
 \href{https://doi.org/10.1215/S0012-7094-04-12723-X}{\textit{Duke Math.~J.}}
 \textbf{127} (2005), 291--311,
 \href{http://arxiv.org/abs/math.QA/0309138}{arXiv:math.QA/0309138}.

\bibitem{MSW11}
Musiker G., Schiffler R., Williams L., Positivity for cluster algebras from
 surfaces, \href{https://doi.org/10.1016/j.aim.2011.04.018}{\textit{Adv.
 Math.}} \textbf{227} (2011), 2241--2308,
 \href{http://arxiv.org/abs/0906.0748}{arXiv:0906.0748}.

\bibitem{MW13}
Musiker G., Williams L., Matrix formulae and skein relations for cluster
 algebras from surfaces,
 \href{https://doi.org/10.1093/imrn/rns118}{\textit{Int. Math. Res. Not.}}
 \textbf{2013} (2013), 2891--2944,
 \href{http://arxiv.org/abs/1108.3382}{arXiv:1108.3382}.

\bibitem{P87}
Penner R.C., The decorated {T}eichm\"uller space of punctured surfaces,
 \href{https://doi.org/10.1007/BF01223515}{\textit{Comm. Math. Phys.}}
 \textbf{113} (1987), 299--339.

\bibitem{S08}
Schiffler R., A cluster expansion formula ({$A_n$} case),
 \href{https://doi.org/10.37236/788}{\textit{Electron.~J. Combin.}}
 \textbf{15} (2008), 64, 9~pages,
 \href{http://arxiv.org/abs/math.RT/0611956}{arXiv:math.RT/0611956}.

\bibitem{Wil19}
Wilson J., Positivity for quasi-cluster algebras,
 \href{http://arxiv.org/abs/1912.12789}{arXiv:1912.12789}.

\end{thebibliography}
\end{document}